\documentclass{article}
\date{}
\usepackage{a4wide}

\usepackage{tikz,pgfplots}
\usetikzlibrary{plotmarks}
\newlength\figureheight
\newlength\figurewidth
\usepackage{subfigure}

\usepackage{amssymb}
\usepackage{amsmath}
\usepackage{mathrsfs}
\usepackage{xcolor}

\usepackage{hyperref}
\hypersetup{
    colorlinks,%
    citecolor=black,%
    filecolor=black,%
    linkcolor=black,%
    urlcolor=black
}

\DeclareMathOperator{\trace}{tr}

\DeclareMathOperator{\rank}{rank}
\usepackage[ruled,vlined]{algorithm2e}

\providecommand{\SetAlgoLined}{\SetLine}

\newcommand{\seqt}[2]{( #1 )_{#2}}
\newcommand{\spa}[1]{\operatorname{span}\left\{#1\right\}}
\newcommand{\spat}[1]{\operatorname{span}\{#1\}}
\newcommand{\inner}[2]{\left\langle #1,#2 \right\rangle}
\newcommand{\innert}[2]{\langle #1,#2 \rangle}
\newcommand{\norm}[1]{\left\|#1\right\|}
\newcommand{\normt}[1]{\|#1\|}
\newcommand{\argmin}[1]{\arg\min_{#1}}

\newcommand{\bN}{\mathbb{N}}

\newcommand{\bR}{\mathbb{R}}

\newcommand{\cC}{\mathcal{C}}

\newcommand{\cH}{\mathcal{H}}
\newcommand{\cI}{\mathcal{I}}

\newcommand{\cL}{\mathcal{L}}
\newcommand{\cM}{\mathcal{M}}

\newcommand{\cT}{\mathcal{T}}
\newcommand{\cU}{\mathcal{U}}
\newcommand{\cV}{\mathcal{V}}
\newcommand{\cW}{\mathcal{W}}
\newcommand{\cX}{\mathcal{X}}

\newcommand{\sP}{\mathscr{P}}

 \title{Low-rank approximate inverse for
   preconditioning tensor-structured linear systems\thanks{This work was
supported by the French National Research Agency (Grant
ANR-2010-COSI-006).}}

\author{L. Giraldi\footnotemark[2]
  \and A. Nouy\footnotemark[2]\ \footnotemark[3] \and G. Legrain\footnotemark[2]}

\begin{document}

\maketitle

\renewcommand{\thefootnote}{\fnsymbol{footnote}}
\footnotetext[2]{Ecole Centrale de Nantes, GeM UMR CNRS 6183, LUNAM Universit\'e, France.}
\footnotetext[3]{Corresponding author (anthony.nouy@ec-nantes.fr).}

\renewcommand{\thefootnote}{\arabic{footnote}}

\begin{abstract}
  In this paper, we propose an algorithm for the construction of
  low-rank approximations of the inverse of an operator given in
  low-rank tensor format. The construction relies on an updated greedy
  algorithm for the minimization of a suitable distance to the inverse
  operator. It provides a sequence of approximations that are defined
  as the projections of the inverse operator in an increasing sequence
  of linear subspaces of operators. These subspaces are obtained by
  the tensorization of bases of operators that are constructed from
  successive rank-one corrections. In order to handle high-order
  tensors, approximate projections are computed in low-rank
  Hierarchical Tucker subsets of the successive subspaces of
  operators. Some desired properties such as symmetry or sparsity can
  be imposed on the approximate inverse operator during the correction
  step, where an optimal rank-one correction is searched as the tensor
  product of operators with the desired properties. Numerical examples
  illustrate the ability of this algorithm to provide efficient
  preconditioners for linear systems in tensor format that improve the
  convergence of iterative solvers and also the quality of the
  resulting low-rank approximations of the solution.
 
\end{abstract}

\pagestyle{myheadings}
\thispagestyle{plain}
\markboth{Loic Giraldi, Anthony Nouy and Gregory Legrain}{Low-rank approximate
  inverse for preconditioning linear systems}

\section{Introduction}\label{sec:intro} 

This paper is concerned with the numerical solution of high-dimensional linear systems of equations in tensor format
\begin{equation}\label{eq:linear_pb}
 A u = b, \quad u\in \bR^{n_{1}}\otimes \hdots \otimes \bR^{n_{d}},
\end{equation}
using low-rank approximation methods. These methods consist in approximating the solution under the form
$$
\sum_{i_{1}} \hdots \sum_{i_{d}}  \alpha_{i_1 \hdots i_{d}} w_{i_{1}}^{1}\otimes \hdots \otimes w_{i_{d}}^{d}, 
$$ 
with $w^{\mu}_{i_{\mu}} \in \bR^{n_\mu}$, $1\le \mu\le d$, and where
the set of coefficients $(\alpha_{i_1 \hdots i_{d}})$ possesses some
particular structure yielding a representation with reduced
complexity. When using suitable approximation formats, low-rank
approximation methods result in a complexity of algorithms that grows
linearly with the dimension $d$, thus allowing the numerical solution
of high-dimensional problems (see the recent surveys
\cite{kolda_tensor_2009,Chinesta2011,Khoromskij2012,Grasedyck2013} and
monograph \cite{Hackbusch2012}). Different strategies have been
proposed for the construction of low-rank approximations of the
solution of equations in tensor format. The first class of methods
consists in defining the approximation as the minimizer in a low-rank
tensor subset of some distance to the solution (e.g. the norm of the
residual of equation \eqref{eq:linear_pb}), see e.g.
\cite{beylkin_algorithms_2005,Espig2012a}. An approximation with
prescribed accuracy can be obtained by introducing an adaptive
selection of tensor subsets or by using greedy constructions where
corrections of the approximation are successively computed in fixed
low-rank subsets (usually rank-one subsets)
\cite{Ammar2010,Cances2011,Falco2011}. A series of improved algorithms
have been proposed in order to increase the quality of suboptimal pure
greedy constructions (see
\cite{Nouy2007,Nouy2008,Nouy2010,Nouy2010-time,Ladeveze2010,Giraldi2012}
and \cite{Falco2012-pgd} for the analysis of a large class of improved
greedy algorithms).
The second class of methods consists in using classical iterative
solvers with low-rank tensor algebra, using efficient algorithms for
low-rank tensor compressions
\cite{kressner_krylov_2010,kressner_low-rank_2011,Ballani2013}.

In this paper, we are interested in the construction of low-rank
preconditioners for equations in tensor format, yielding
preconditioned equations 
 $$PA u =Pb$$ with a preserved low-rank tensor format. 
 Preconditoning aims at improving  the convergence of iterative
 methods but also at improving the quality of low-rank approximations
 defined from the residual of the equation. Different strategies have
 been proposed for the construction of low-rank preconditioners.  In
 the case of equations resulting from a discretization of stochastic
 equations, a rank-one preconditioner has been introduced  in
 \cite{ghanem_numerical_1996}. It is based on the inverse of the
 expectation of the random operator, and it is particularly efficient
 when the random operator has a small variance. In
 \cite{langville_kronecker_2004}, a more general rank-one
 preconditioner has been defined as the inverse of a rank-one
 approximation of the operator. This preconditioner has 
 been exploited in \cite{ullmann_kronecker_2010} for the solution of
 equations arising from the discretization of stochastic parametric
 equations. In the same context, a rank-one preconditioner has
 also been defined in \cite{zander_tensor_2012} as the solution of the
 minimization of $\normt{I-PA}$ over the set of rank-one operators $P$.

 Rank-one preconditioners may be efficient if the operator $A$ only
 slightly deviates from a rank-one operator. In order to address more
 general situations, different strategies have been proposed for the
 construction of higher rank preconditioners. In
 \cite{khoromskij_tensor-structured_2009}, a preconditioner is
 obtained by truncating an expansion of the inverse of the operator.
 In \cite{touzene_tensor_2008}, a preconditioner $P$ is defined as the
 best approximation of the inverse of the operator with the particular
 structure $P = P^1 \otimes I \otimes \hdots \otimes I + \hdots + I
 \otimes \hdots \otimes I \otimes P^d$ corresponding to a rank-$d$
 preconditioner. More recently, an algorithm has been proposed in
 \cite{oseledets_solution_2012} for the construction of a low-rank
 preconditioner $P$ in tensor-train format. It relies on the solution
 of the equation $AP=I$ with a DMRG algorithm, this algorithm allowing
 for an automatic selection of the rank. In order to avoid the
 inversion of large matrices (large $n_{\mu}$), a quantization
 technique is introduced.

 In the present paper, we propose an algorithm for the computation of a
 low-rank approximation $P$ of the inverse operator $A^{-1}$ using
Tucker or Hierarchical Tucker format. This algorithm is an updated greedy
 algorithm for the minimization of a suitable distance $\Vert
 A^{-1}-P\Vert_{\star}$. The norm $\Vert\cdot\Vert_{\star}$ is chosen such
 that the approximation can be computed without any a priori
 approximation of $A^{-1}$, and it is chosen according to the
 properties of $A$ (namely symmetric positive definite or simply
 definite operator). Compared to a direct minimization of $\Vert
 A^{-1}-P\Vert_{\star}$ over a set of Tucker or Hierarchical Tucker tensors with
 given rank, the greedy procedure has the advantages of being adaptive
 and of considerably reducing the complexity of the construction of a
 low-rank approximation, therefore allowing the manipulation of large
 dimensions $n_{\mu}$. Starting from $P_{0}=0$, one step of the
 updated greedy algorithm consists in (i) computing a rank-one
 correction of the previously computed approximation $P_{r-1}$ by
 minimizing $\Vert A^{-1}-P_{r-1}-W_{r}\Vert_{\star}$ over the set of
 rank-one operators $W_{r}=W^{1}_{r}\otimes \hdots \otimes W^{d}_{r}$,
 (ii) updating reduced spaces of operators $\cU^{\mu}_{r}$ ($1\le
 \mu\le d$) which are defined as the span of the set of operators
 $\{W^{\mu}_{1},\hdots,W^{\mu}_{r}\}$, and (iii) computing a new
 approximation $P_{r}$ in the space $\cU_{r} = \cU^{1}_{r}\otimes
 \hdots \otimes \cU^{d}_{r}$ by minimizing $\Vert
 A^{-1}-P_{r}\Vert_{\star}$ in $\cU_{r}$ or over a set of low-rank Hierarchical Tucker
 tensors in $\cU_{r}$. More precisely, the approximation $P_{r}$ is
 searched under the form
 $$P_{r} = \sum_{i_{1}=1}^{r}\hdots \sum_{i_{d}=1}^{r}
 \alpha_{i_{1},\hdots,i_{d}} W_{i_{1}}^{1}\otimes \hdots \otimes
 W_{i_{d}}^{d},$$ where the set of coefficients $\alpha$ is optimized in $\mathbb{R}^{r}\otimes
 \hdots \otimes \mathbb{R}^{r}$ or in a low-rank Hierarchical Tucker subset of $\mathbb{R}^{r}\otimes
 \hdots \otimes \mathbb{R}^{r}$. For the solution of the minimization
 problems over the set of rank-one tensors (step (i)) and the set of
 Hierarchical Tucker tensors with bounded rank (step (iii)),
 alternating minimization algorithms are used \cite{kressner_preconditioned_2011,Uschmajew2013}. 
 \\
 Some desired properties such as symmetry and sparsity can be imposed
 on the approximate inverse. This is done in the correction
 step (i) where the optimal rank-one correction is searched as
 the tensor product of operators with the desired properties. During
 the alternating minimization algorithm, imposing the symmetry on the
 matrix $W^{\mu}_{r}$ requires the solution of a Sylvester equation. For imposing sparsity on $W^{\mu}_{r}$, we propose a
 straightforward generalization of the sparse approximate inverse
 algorithm proposed in \cite{grote_parallel_1997}, which is an
 adaptive algorithm for the determination of the sparsity pattern.

\par
The outline of the paper is as follows. In Section \ref{sec:tensor},
we briefly recall some useful definitions on tensor spaces and
low-rank tensor approximations. In Section \ref{sec:rank-one-approx}, we introduce an algorithm for computing a rank-one approximation of the inverse operator, with possible imposed properties.  
 In Section
\ref{sec:updated_greedy_algorithm}, we introduce the algorithm for computing a low-rank approximate inverse in low-rank Tucker or Hierarchical Tucker formats. 
In Section \ref{sec:applications}, the efficiency of the proposed
preconditioning technique is illustrated on numerical problems: a
Poisson equation in high dimension (symmetric problem) and a
 linear equation resulting from the discretization of a
stochastic partial differential equation using spectral stochastic
methods.

\section{Tensor spaces and low-rank tensor approximation} \label{sec:tensor} 

\subsection{Tensor spaces}
Let $D = \{1,\hdots,d\}$, with 
 $d \in \bN^*$. Let $\cX^{\mu}$, $\mu \in D$, be a finite dimensional space equipped with an inner product 
  $\innert{\cdot}{\cdot}_\mu$ and the
associated norm $\normt{\cdot}_\mu$.  We consider the tensor space $\cX= \cX^{1}\otimes \hdots \otimes \cX^{d}$. 
In the following, the simplified notation  $\bigotimes_\mu$ will be
used for $\bigotimes_{\mu \in D}$, as well as $\bigotimes_{\mu \ne
  \lambda}$ for $\bigotimes_{\mu \in D \setminus \{\lambda\}}$, $\lambda\in D$.
A tensor $x \in \cX$ can be written under the form
$x = \sum_{i=1}^r \bigotimes_\mu x_i^\mu$ for some  $ r \in
\bN$ and $\ x_i^\mu \in \cX^{\mu}$.
The minimal integer $r$ which allows to represent $x$ exactly under this form is called the (canonical)
rank of $x$. $\cX$ is a Hilbert space for the induced inner product
$\innert{\cdot}{\cdot}$ defined for rank-one tensors by 
$\langle{\bigotimes_\mu x^\mu} ,  {\bigotimes_\mu y^\mu} \rangle = \prod_{\mu \in D}
\langle{x^\mu},{y^\mu} \rangle_\mu$,
and extended by linearity to the whole space $\cX$. The associated norm is noted 
$\normt{\cdot}$.

For $\cX^{\mu}=\bR^{n_{\mu}}$, we use the vector 2-norm $\normt{\cdot}_\mu$ and the associated 
canonical inner product  $\innert{\cdot}{\cdot}_\mu$. For $\cX^{\mu}=\bR^{n_{\mu}\times n_{\mu}}$, we take for  $\normt{\cdot}_\mu$ the Frobenius norm and for $\innert{\cdot}{\cdot}_\mu$ the associated inner product defined by 
$ \innert{A^{\mu}}{B^{\mu}}_\mu =
\trace((A^{\mu})^{T} B^{\mu})$ for $A^{\mu},B^{\mu}\in \bR^{n_\mu \times
  n_\mu}$. 
A matrix $A^{\mu}\in\bR^{n_\mu \times
  n_\mu} $ is identified with the corresponding operator $A^{\mu} : \bR^{n_{\mu}} \rightarrow \bR^{n_{\mu}}$, and a tensor  $A\in \bigotimes_{\mu} \bR^{n_\mu \times
  n_\mu}  $ is identified with the corresponding operator $A : \bigotimes_{\mu} \bR^{n_\mu} \rightarrow 
  \bigotimes_{\mu} \bR^{n_\mu}$. We denote by $A^{T}$ the adjoint of $A$ for the induced inner product, which for $A=\sum_{i=1}^r \bigotimes_\mu A_i^\mu$ is obtained by $A^{T} =\sum_{i=1}^r \bigotimes_\mu (A_i^\mu)^{T} $, where $(A_i^\mu)^{T} $ denotes the  transpose of matrix $A_{i}^{\mu}$.

\subsection{Low-rank tensor approximation}

Let $\cM\subset \cX$ be a subset of low-rank tensors in $\cX$. The best approximation of $y$  in $\cM$, if it exists, is defined by 
\begin{align*}
\min_{x \in \cM} \norm{y- x}.
\end{align*}
We also define a quasi-best approximation of $y$ in $\cM$ as an element $x^{\gamma}\in \cM$ satisfying 
\begin{align*}
  \norm{y- x^{\gamma}} \le \gamma \inf_{x \in \cM} \norm{y - x}
\end{align*}
for some factor $\gamma>1$. Note that a quasi-best approximation exists for any $\gamma>1$. 

\subsubsection{Canonical format}
The subset of rank-$r$ canonical tensors is defined by 
\begin{equation*}
  \cC_r(\cX) = \left\{ x = \sum_{i=1}^r \bigotimes_\mu x_i^\mu;\ x_i^\mu \in
\cX^{\mu},\ \forall \mu \in D \right\}.
\end{equation*}
{This format is simple. However, the set $\cC_r(\cX)$ is not closed for
$d>2$ and $r>1$, therefore the best approximation problem using
canonical format is ill-posed in this case \cite{de_silva_tensor_2008}}. Alternating minimization algorithm \cite{kolda_tensor_2009} or other optimization algorithms \cite{acar_scalable_2011,espig_regularized_????} can be used to solve the approximation problem in $\cC_r(\cV)$. However, for $d>2$, there is no available algorithm for computing a quasi-best approximation 
in $\cC_r(\cX)$ with a controlled factor $\gamma$. 
 
 \subsubsection{Tucker format}

The subset of Tucker tensors with rank $r=(r_{\mu})_{\mu\in D}$, introduced in 
\cite{tucker_mathematical_1966}, is defined by 
\begin{equation}\label{eq:tucker_set}
  \cT_r(\cX) = \left\{x \in \cX ;
\begin{array}{l} \displaystyle \text{there exist linear subspaces $\cU^\mu$ with } \\
\text{$dim(\cU^\mu)=r_\mu$,  $1\le \mu\le d$, such that }x\in
\bigotimes_{\mu=1}^d \cU^{\mu}
  \end{array}
 \right\}.
\end{equation}
Letting $\cI_{r} = \cI_{r_{1}}\times\hdots \times \cI_{r_{d}}$ with $\cI_{r_{\mu}} = \{1,\hdots,r_{\mu}\}$, it is equivalently defined by
\begin{equation}\label{eq:tucker_set_{bis}}
  \cT_r(\cX) = \left\{x = \sum_{i \in \cI_{r}} \alpha_i \bigotimes_\mu x_{i_\mu}^\mu \; ;\;
  \alpha \in \bigotimes_{\mu=1}^{d} \bR^{r_{\mu}}, \;
x_{i_\mu}^\mu \in \cX^{\mu} 
 \right\},
\end{equation}
where $\alpha$ is the so called core tensor of the Tucker representation. 
The set $\cT_r(\cX)$ is closed \cite{falco_minimal_2010}, so that a
best approximation problem of a tensor in $ \cT_r(\cX)$ always exists.
Moreover, when using the canonical norm, the Higher Order Singular
Value Decomposition (HOSVD) algorithm proposed in
\cite{de_lathauwer_multilinear_2000} allows the efficient computation
of a quasi-best approximation of a tensor with a controlled factor
$\gamma = \sqrt{d}$. This quasi-best approximation can further be
improved using Higher Order Orthogonal Iterations (HOOI) algorithm
\cite{DeLathauwer2000}, which is an alternating minimization
algorithm. The drawback of this format is that the core tensor
$\alpha$ is of order $d$ so that this format suffers from the curse of
dimensionality.

\subsubsection{Hierarchical Tucker format}\label{sec:hierarchical-tensors}

The Hierarchical Tucker tensor format has been introduced in
\cite{hackbusch_new_2009} and is defined as follows. Let $T$ be a tree
of dimensions, defined as a full binary tree on $D$ with root $D$. The
set of leaves of the tree is defined by $L(T) 
=\{\{\mu\};\mu\in D\}$ and the set of interior nodes is $I(T) = T
\setminus L(T)$. The set of successors $S(t)$ of an interior note
$t\in I(T)$ is composed by two nonempty successors $t_1$ and $t_2$ in
$T$ such that $t = t_1 \cup t_2$ and $t_1 \cap t_2 = \emptyset$. The
complement of $t \in T$ is denoted by $t^c = D\setminus t$. We denote
by $\cX^t = \bigotimes_{\mu \in t} \cX^{\mu}$ and $\cX^{t^c} =
\bigotimes_{\mu \in t^c} \cX^{\mu}$. The $t$-matricization operator
$M_t : \cX \to \cX^t \otimes \cX^{t^c}$ is defined for $ x =
\sum_{i=1}^{r} \bigotimes_\mu x_i^\mu$ by
\begin{equation*}
M_t(x) = \sum_{i=1}^r
\left(\bigotimes_{\mu \in t} x_i^\mu\right) \otimes
\left(\bigotimes_{\mu \in t^c} x_i^\mu\right).
\end{equation*}
The $t$-rank of $x$ is then defined as the rank of the 2-order tensor 
$M_{t}(x)\in \cX^t \otimes
\cX^{t^c}$.
Given a tree of dimensions $T$ and a family of ranks $r = \seqt{r_t}{t \in
  T}$ associated with the tree, the set of Hierarchical Tucker tensors
with bounded rank $r$ is defined by
\begin{equation*}
  \cH_r^T(\cX) = \left\{ x \in \cX;\ \text{$t$-$\rank$}(x) \le r_t,\ \forall t
    \in T\right\}.
\end{equation*}
We note that the set of Tucker tensors $\cT_{(r_{1},\hdots,r_{d})}= \left\{ x \in \cX; \
    \text{$\mu$-$\rank$}(x) \le r_{\mu},\ \forall \mu \in D \right\}$, and therefore
    $$\cH_r^T(\cX) \subset \cT_{(r_{1},\hdots,r_{d})}(\cX).$$
In practice, a tensor $x \in \cH_r^T(\cX)$ can be represented
 under the form
\begin{align*}
  x &= \sum_{i=1}^{r_{t_1}} \sum_{j=1}^{r_{t_2}} \beta_{ij}^D x_i^{t_1}
  \otimes x_j^{t_2},\; \beta^D \in \bR^{r_{t_1}\times r_{t_2}},\; \{t_1, t_2\}=S(D),
  \end{align*}
where for all $t\in I(T) \setminus\{D\}$, 
\begin{align*}
  x_k^t &= \sum_{i=1}^{r_{t_1}} \sum_{j=1}^{r_{t_2}} \beta_{ijk}^t x_i^{t_1}
  \otimes x_j^{t_2}, \; \beta^t  \in \bR^{r_{t_1}\times r_{t_2}\times r_t}, \; 
  \{t_1,t_2\}=S(t), \; k\in\{1,\hdots,r_{t}\}.\end{align*}
Therefore, $x \in \cH_r^T(\cX)$ is completely determined by 
 the set of transfer tensors $\{\beta^t ; t \in I(T)\}$ associated with the interior nodes of the tree, and the set of elements $\{x_{i}^\mu ; i\in \cI_{r_{\mu}} ,  \mu\in D\}$ associated with the leaves of the tree. The tensor $x$ can be written 
  $$
 x = \sum_{i\in \cI_{r_{1}}\times \hdots \times \cI_{r_{d}}} \alpha_{i} \bigotimes_{\mu=1}^{d} x_{i_{\mu}}^{\mu},  $$
 where the tensor $\alpha \in \cH_{r}^{T}(\bigotimes_{\mu=1}^{d} \bR^{r_{\mu}})$ is a rank-$r$ Hierarchical Tucker tensor that can be expressed in terms of the 
 transfer tensors $\{\beta^t ; t \in I(T)\}$. 
 
 The set $\cH_r^T(\cX)$ is closed (see \cite{Hackbusch2012},
section 11.4.1.1) so that a best approximation of a tensor in
$\cH_r^T(\cX)$ always exists. Moreover, when using the induced canonical norm, the Hierarchical Singular Value Decomposition (HSVD) algorithm proposed in \cite{grasedyck_hierarchical_2010} allows the  efficient computation of a quasi-best approximation of a tensor in $\cH_r^T(\cX)$ with a controlled factor $\gamma = \sqrt{2d-3}$.
Besides, given that $\#I(T) = d-1$
and $\#L(T) = d$, the Hierarchical Tucker format does not suffer from the curse of
dimensionality since the dimension  of the parametrization of a tensor in $\cH_r^T(\cX)$ depends linearly on $d$.

\section{Rank-one approximation of an inverse operator} 
\label{sec:rank-one-approx}

Let consider $A $ in $ \cW = \bigotimes_{\mu=1}^{d} \cW^{\mu}$, with $\cW^{\mu} = \bR^{n_{\mu}\times n_{\mu}}$.
In this section, we introduce an algorithm for computing a rank-one approximation of the inverse $A^{-1}$ of $A$, considered as an operator from $\bigotimes_{\mu=1}^{d}\bR^{n_{\mu}}$ to $\bigotimes_{\mu=1}^{d} \bR^{n_{\mu}}$. The rank-one approximation is 
searched in a linear subspace $\cU = \bigotimes_{\mu=1}^{d} \cU^{\mu}$ of $\cW$, with  $\cU^{\mu} =\cW^{\mu}$ or $\cU^{\mu}\subset \cW^{\mu}$, a linear subspace of operators  with prescribed properties such as symmetry or sparsity. 

 \subsection{Best rank-one correction}
 
Let $P \in \cW$ be a first approximation of $A^{-1}$ (e.g. a known preconditioner of $A$).
$\cC_{1}(\cU)$ denotes the set of rank-one elements of the tensor space  $\cU \subset \cW$. 
The best rank-one correction $W = \bigotimes_{\mu=1}^{d} W^{\mu}\in \cC_{1}(\cU)$ of $P$ is defined by the following problem 
\begin{equation}
\min_{W\in \cC_{1}(\cU)} \Vert A^{-1} - (P + W)\Vert_{\star}\label{eq:min-rank-one-star}
\end{equation}
where $\Vert\cdot\Vert_{\star}$ is a norm on $\cW$ that has to be chosen such that an approximation $W$ can be computed without knowing $A^{-1}$. 
 
 \subsection{Definition of the norm $\Vert \cdot\Vert_{\star}$}\label{sec:norm-star}
 
If $A$ is \emph{symmetric positive definite}, we choose the norm $\Vert\cdot\Vert_{\star}=\Vert\cdot\Vert_{A}$ defined by $ 
\Vert X \Vert_{A} = \Vert XA^{1/2}\Vert = \sqrt{\innert{XA}{X}}$ and associated with the inner product 
$\innert{\cdot}{\cdot}_{A}$ defined by
$\innert{X}{Y}_A = \innert{XA}{Y}$. We note that $\Vert  A^{-1} - (P+W) \Vert_{A} = \Vert I - (P+W)A \Vert_{A^{-1}},$ so that the minimization 
problem \eqref{eq:min-rank-one-star}
provides  	a left approximate inverse $(P+W)$ of $A$.\footnote{We could also choose $\Vert \cdot \Vert_{\star} = \Vert \cdot \Vert_{A} $ with $\Vert X \Vert_{A} =  \Vert A^{1/2} X\Vert = \sqrt{\innert{AX}{X}}$, so that the minimization of $\Vert  A^{-1} - (P+W) \Vert_{\star} = \Vert I - A(P+W) \Vert_{A^{-1}}$ provides a right approximate inverse  $(P+W)$  of $A$.} 

In the more general case of a \emph{definite operator} $A$,
we choose the norm $\Vert\cdot\Vert_{\star}=\Vert\cdot\Vert_{AA^{T}}$ defined by $ 
\Vert X \Vert_{AA^{T}} = \Vert XA \Vert = \sqrt{\innert{XA}{XA}} = \sqrt{\innert{XAA^{T}}{X}} $ and associated with the inner product 
$\innert{\cdot}{\cdot}_{AA^{T}}$ defined by
$\innert{X}{Y}_{AA^{T}} = \innert{XA}{YA} =  \innert{XAA^{T}}{Y}$. 
We note that $\Vert  A^{-1} - (P+W) \Vert_{AA^{T}} = \Vert I - (P+W)A \Vert,$ so that the minimization 
problem \eqref{eq:min-rank-one-star}
provides a left approximate inverse $(P+W)$ of $A$.\footnote{We could also choose $\Vert X \Vert_{\star} = \Vert A X\Vert = \sqrt{\innert{AX}{AX}}$, so that the minimization of $\Vert  A^{-1} - (P+W) \Vert_{\star} = \Vert I - A(P+W) \Vert$ provides a right approximate inverse  $P+W$  of $A$.}

From now on, we consider that $\Vert \cdot\Vert_{\star} = \Vert \cdot \Vert_{AB}$, with 
$B=I$ if $A$ is symmetric positive definite or
$B=A^{T}$ if $A$ is simply definite,
and we denote by  $\innert{\cdot}{\cdot}_{\star} = \innert{\cdot}{\cdot}_{AB}$ the associated inner product on $\cW$. 

\subsection{Stationarity conditions}

A necessary condition of optimality for a solution $W=\bigotimes_{\mu=1}^{d }W^{\mu}$  of problem  \eqref{eq:min-rank-one-star} is 
\begin{equation}\label{eq:stationarity-equation-star}
\inner{A^{-1}-(P + W)}{\delta W}_{\star}=
0, 
\quad \forall \delta W \in T_W(\cC_1(\cU)),
\end{equation}
where $T_W(\cC_1(\cU))$ is the tangent space of $\cC_{1}(\cU)$ at $W$ defined by
 \begin{align*}
& T_W(\cC_1(\cU)) = \sum_{\mu =1}^{d} T_W^{(\mu)}(\cC_1(\cU)) \subset \cC_{d}(\cU) ,
\end{align*}
with
 \begin{align*} 
 T_W^{(\mu)}(\cC_1(\cU)) = 
\left\{ \delta W = W^{1} \otimes \hdots \otimes 
 \delta W^{\mu}
 \otimes \hdots \otimes W^{d} \in \cC_{1}(\cU) ; 
 \delta W^{\mu} \in \cU^{\mu}
\right\}.
\end{align*}
Using the definition of inner  product $\innert{\cdot}{\cdot}_{\star}=\innert{\cdot}{\cdot}_{AB}$, we have that $W$ satisfies \eqref{eq:stationarity-equation-star} if and only if 
\begin{equation}\label{eq:stationarity-equation}
\inner{B-(P + W)AB}{\delta W} = 
0, 
\quad \forall \delta W \in T_W(\cC_1(\cU)).
\end{equation}
Let $\sP^{\mu} : \cW^{\mu} \rightarrow \cU^{\mu}$, $\mu\in D$, denote the orthogonal projector from $\cW^{\mu}$ into $\cU^{\mu}$, such that $\sP^{\mu}\sP^{\mu}=\sP^{\mu}$ and 
$
\innert{\sP^{\mu}(X^{\mu})}{Y^{\mu}-\sP^{\mu}(Y^{\mu})}_{\mu} = 0$ for all  $X^{\mu},Y^{\mu}\in \cW^{\mu}
$. 
The operator $\sP :\cW \rightarrow \cU$ defined by $\sP = \bigotimes_{\mu=1}^{d} \sP^{\mu}$ is the orthogonal projector from 
$\cW$ to $\cU$ such that $\sP\sP = \sP$ and
\begin{equation}\label{projector-ortho}
\inner{\sP(X)}{Y-\sP(Y)} = 0,\quad \forall X,Y\in \cW.
\end{equation}
Noting that $T_W(\cC_1(\cU))\subset \cU$, we have that for all $\delta W\in T_W(\cC_1(\cU))$, 
\begin{align*}
  \inner{B - (P+W)AB}{\delta W} &= \inner{B - (P+W)AB}{\sP(\delta W)}\\
 &= \inner{\sP(B - (P+W)AB)}{\sP(\delta W)} \\
 &= 
 \inner{\sP(B - (P+W)AB)}{\delta W}.
 \end{align*}
Therefore, using the definition of the tangent space $T_W(\cC_1(\cU))$, we obtain that $W$ satisfies \eqref{eq:stationarity-equation} if and only if for all $\lambda\in D$,
\begin{equation}
  \label{eq:PWAB}
  \inner{\sP(WAB)}{\delta W} = \inner{\sP(R(P))}{\delta
    W}, \quad \forall \delta W \in T_W^{(\lambda)}(\cC_1(\cU)),
\end{equation}
with $R(P) = B-PAB$.

\subsection{Alternating minimization algorithm}\label{sec:altern-minim}
For solving \eqref{eq:min-rank-one-star}, we use an alternating minimization algorithm. Starting 
from an arbitrary initialization $W = \otimes_{\mu} W^{\mu} \in \cC_{1}(\cU)$, it consists in solving 
successively   the quadratic optimization problems
\begin{equation}\label{eq:ama-min-mu}
\min_{W^{\lambda} \in\cU^{\lambda}}  \Vert A^{-1}-P-\bigotimes_{\mu=1}^{d} W^{\mu} \Vert_{\star}^{2},
 \end{equation}
 for $\lambda \in \{1,\hdots,d,1,\hdots \}$. 
It is observed that this algorithm converges to an element  $W$ that satisfies the stationarity condition \eqref{eq:stationarity-equation-star}, or equivalently \eqref{eq:PWAB} for all $\lambda\in D$. 
 However, note that it does not necessarily yield a solution of \eqref{eq:min-rank-one-star}.

The solution of \eqref{eq:ama-min-mu} is equivalent to the solution of Equation  \eqref{eq:PWAB} for given $W^{\mu}$, $\mu\neq \lambda$.
Therefore, noting $ C = AB = \sum_{i=1}^{r_C} \otimes_\mu C_i^\mu$ and 
 $R(P) = B-PAB =  \sum_{i\in \cI_{r_{R}}} \gamma_{i} \otimes_{\mu} R_{i_{\mu}}^{\mu} $,  we obtain that \eqref{eq:ama-min-mu}  is 
 equivalent to the following linear problem 
 on $W^{\lambda}\in \cW^{\lambda}$:  
 \begin{align}\label{eq:proj-lambda}
  \sP^\lambda(W^\lambda Q^\lambda) = \sP^\lambda(H^\lambda(P)),
\end{align}
with
\begin{align}
  &Q^\lambda =  \sum_{i=1}^{r_C} C_i^\lambda \prod_{\mu\neq \lambda}
  \inner{\sP^\mu(W^\mu C_i^\mu)}{W^\mu}_\mu   , \label{eq:proj-lambda-Q} \\
  &H^\lambda(P)= \sum_{i\in \cI_{r_{R}}} R_{i_{\lambda}}^\lambda \gamma_{i} \prod_{\mu\neq \lambda}
  \inner{\sP^\mu(R_{i_{\mu}}^\mu)}{W^\mu}_\mu .\label{eq:proj-lambda-H}
\end{align}

If $\cU^{\lambda} = \cW^{\lambda}$, that means if we do not impose any
particular property to matrices associated with dimension $\lambda$, then $\sP^{\lambda} $ is 
the identity on $\cW^{\lambda}$ and Equation  \eqref{eq:proj-lambda} becomes
\begin{align}\label{eq:proj-lambda-no-constraint}
W^\lambda Q^\lambda = H^\lambda(P),
\end{align}
with $Q^{\lambda}$ and $H^{\lambda}(P)$ defined in Equations 
\eqref{eq:proj-lambda-Q} and \eqref{eq:proj-lambda-H} respectively.

\subsection{Imposing properties}
\label{sec:imposing-properties}

\subsubsection{Imposing symmetries}
\label{sec:imposing-symmetries}

We note $D_{sym}$, $D_{skew}$ and $D_c$ three sets of indices
such that they form a partition of $D = \{1,\hdots,d\}$.
Then, we consider $\cU = \bigotimes_\mu
\cU^\mu$ such that
\begin{equation*}
  \cU^\mu =
  \begin{cases}
    \{X \in \cW^{\mu};\ X = X^T\} &\mbox{if } \mu \in D_{sym}, \\
    \{X \in  \cW^{\mu};\ X = -X^T\} &\mbox{if } \mu \in D_{skew}, \\
     \cW^{\mu} &\mbox{if } \mu \in D_c.
  \end{cases}
\end{equation*}
The orthogonal projector $\sP^\mu: \cW^{\mu} \rightarrow \cU^\mu$ is such that 
\begin{equation*}
  \sP^\mu(X) =
  \begin{cases}
    \frac{1}{2} (X+X^T) &\mbox{if } \mu \in D_{sym}, \\
    \frac{1}{2} (X-X^T) &\mbox{if } \mu \in D_{skew}, \\
    X &\mbox{if } \mu \in D_c.
  \end{cases}
\end{equation*}
Therefore, the linear equation \eqref{eq:proj-lambda} in $W^{\lambda}\in \cU^{\lambda}$ 
can  be written
\begin{align*}
  \begin{cases} W^\lambda Q^\lambda + \left(Q^\lambda\right)^T W^\lambda =
  H^\lambda + \left(H^\lambda\right)^T&\mbox{if } \lambda \in D_{sym},\\
  W^\lambda Q^\lambda - \left(Q^\lambda\right)^T W^\lambda =
  H^\lambda - \left(H^\lambda\right)^T&\mbox{if } \lambda \in D_{skew}, \\  
   W^\lambda Q^\lambda = H^\lambda, &\mbox{otherwise}.
\end{cases}
\end{align*}
We notice that the equations associated with $\lambda$ in $ D_{sym}$ or 
$D_{skew}$ are particular cases of the so called Sylvester equation
which can be solved with the algorithm from
\cite{bartels_algorithm_1972}. 

If $A$ and $P$ are symmetric, and if $D_{sym} = D$ (that means that we search for a 
symmetric rank-one approximation), 
then $Q^\lambda$ is  symmetric and  \eqref{eq:proj-lambda} is a 
continuous Lyapunov equation
\begin{align*}
  W^\lambda Q^\lambda + Q^\lambda W^\lambda = H^\lambda + (H^\lambda)^{T}.
\end{align*}

\subsubsection{Imposing sparsity}
\label{sec:imposing-sparsity}

Here, we are interested in using sparse approximation in order to handle large matrices. 
For $\lambda\in D$, let $I^{\lambda} \subset \{1,\hdots,n_{\lambda}\}^{2}$ and let $\cU^{\lambda}$ be the 
subspace of matrices with sparsity pattern 
$I^{\lambda}$:
\begin{align*}
\cU^\lambda= \left\{ X\in \cW^{\lambda} ; (X)_{kj} = 0,\ \forall
(k,j) \notin I^\lambda\right\}.
\end{align*}
The orthogonal projector from $\cW^{\lambda}$ onto $\cU^{\lambda}$  is 
defined for $X^\lambda\in \cW^{\lambda}$ by
\begin{align*}
(\sP^\lambda(X^\lambda))_{kj}=\begin{cases} (X^\lambda)_{kj}   & \text{if }
(k,j) \in I^\lambda, \\
0 & \text{if } (k,j) \notin I^\lambda. \\
\end{cases}
\end{align*}
Similarly to the SParse Approximate Inverse (SPAI) method from  \cite{grote_parallel_1997}, we reformulate equation \eqref{eq:proj-lambda} 
as the following minimization
problem:
\begin{align}\label{eq:funct_spai}
  \min_{W^\lambda\in \cU^{\lambda}} \norm{\sP^\lambda(W^\lambda Q^\lambda -
    H^\lambda)}^2_\lambda.
\end{align}
Noting $\{{w_k^\lambda}\}_{1 \le k \le n_\lambda}$
(resp. $\{{h_k^\lambda}\}_{1 \le k \le n_\lambda}$)  the rows of
$W^\lambda$ (resp. $H^\lambda$), we have 
\begin{align*}
  \norm{\sP^\lambda(W^\lambda Q^\lambda-H^\lambda)}_\lambda^2 =
  \sum_{k=1}^{n_\lambda} \norm{w_k^\lambda  Q^\lambda P_k^\lambda-
    h_k^\lambda P_k^\lambda}_\lambda^2,
\end{align*}
where $P_{k}^{\lambda}\in  \bR^{n_{\lambda}\times n_{\lambda}}$ is a  boolean diagonal matrix such that $(P_{k}^{\lambda})_{jj} =1$ if $j\in I^{\lambda}_{k} $ and $(P_{k}^{\lambda})_{jj} =0$ if $j\notin I^{\lambda}_{k}$, where $I^{\lambda}_{k} = \{ j ; (k,j) \in I^{\lambda}\}$ denotes the pattern of the row $k$.  
Therefore, the minimization problem  \eqref{eq:funct_spai} is equivalent to $n_{\lambda}$
independent minimization problems  (that can be solved in parallel)
\begin{align}\min_{w_{k}^{\lambda}\in\bR^{n_{\lambda}}}\norm{w_k^\lambda
   Q^\lambda P_k^\lambda- h_k^\lambda P_k^\lambda}^2_\lambda \quad \text{submitted to } 
  (w_{k}^{\lambda})_{j} =0 \; \text{ for all } \; j \notin I^{\lambda}_{k}.\label{eq:spai-probk}
 \end{align}
Each problem \eqref{eq:spai-probk} can be rewritten 
as the minimization of $ \Vert{\widehat w_k^\lambda \widehat
  Q^\lambda_{k} - \widehat h_k^\lambda}\Vert_\lambda^2$ over $\widehat w_{k}^{\lambda}\in \bR^{m^{\lambda}_{k}}$,  
  with $m^{\lambda}_{k} = \# I^{\lambda}_{k}$, where $\widehat w_{k}^{\lambda}$ (resp. $\widehat h_{k}^{\lambda}$) denotes the vector of non zero entries of $w^{\lambda}_{k}$ (resp. $h^{\lambda}_{k}P_k^\lambda$).  This minimization problem on $\widehat w_{k}^{\lambda}$ can be solved
using a QR decomposition of the reduced matrix $\widehat Q^\lambda_{k}$ (see  \cite{grote_parallel_1997}). 
\\
For the adaptive selection of the pattern $I^{\lambda}$, we use the
iterative method proposed in  \cite{grote_parallel_1997}. Let $\sP^{\lambda,(i)}$ be the projector associated to the set of patterns $\{I_k^{\lambda,(i)}\}_{1 \le k \le n_{\lambda}}$. We start 
from an initial projector $\sP^{\lambda,(0)}$ (e.g. associated with diagonal patterns $I_k^{\lambda,(0)} = \{k\}$). Then, for $i=0,\hdots,i_{\max}$, we proceed as follows. We compute
\begin{equation*}
  W^{\lambda,(i)} = \arg \min_{W^\lambda\in \cU^{\lambda}} \norm{\sP^{\lambda,(i)}(W^\lambda Q^\lambda -
    H^\lambda)}^2_\lambda.
\end{equation*}
Then, for each row $k$, we compute a new index 
\begin{align*}
 j^{\lambda,(i)}_{k} = \arg \min_{1 \le j \le n_\lambda} \left( \min_{\gamma \in \bR}
  \norm{(w_k^{\lambda,(i)} + \gamma (e_j^\lambda)^{T})  Q^\lambda - h_k^\lambda }_\lambda^2\right),
\end{align*}
where $e_j^{\lambda}$ is the $j$-th canonical basis vector in $\bR^{n_{\lambda}}$, and we set $I_k^{\lambda,(i+1)} = I_k^{\lambda,(i)} \cup
\{j^{\lambda,(i)}_{k}\}$.

\section{A constructive algorithm using projections in reduced spaces}\label{sec:updated_greedy_algorithm}
Here, we introduce a constructive algorithm for the approximation of the inverse $A^{-1}$ of an operator $A$ using Tucker or Hierarchical Tucker tensor format. 

The algorithm starts with an initialization $P_0$, such as $P_{0}=0$ or a first approximation of $A^{-1}$. 
Knowing an approximation $P_{r-1}$ of
$A^{-1}$, we look for a rank-one correction $W_{r}=\bigotimes_\mu W_r^\mu$ which solves
\begin{align}
  \label{eq:greedy}
  \min_{W \in \cC_1(\cW)}\norm{A^{-1} - P_{r-1}-W}_{\star},
\end{align}
where the norm $\Vert\cdot\Vert_{{\star}}$ has been defined in Section \ref{sec:norm-star}.
 Then, we define the linear subspace 
 \begin{equation*}
   \cU_r = \bigotimes_{\mu=1}^{d} \cU_r^\mu  \qquad \text{with} \qquad\cU_r^\mu = \spa{W_1^\mu, \hdots, W_r^\mu}.
\end{equation*}
The dimension of $\cU_{r}$ is $\prod_{\mu} r_{\mu}$, with $r_{\mu}\le r$.
In practice, we construct orthogonal bases $\{Q_i^\mu\}_{1 \le i \le r_\mu}$ of the subspaces
$\cU_r^\mu$, yielding a basis $\{ \bigotimes_{\mu}Q_i^\mu\}_{i\in \cI_{r}}$ of $\cU_{r}$, with $\cI_{r} = \left\{(i_1,\hdots,i_d)\in \bN^d;
  \ 1 \le i_\mu \le r_\mu \right\}$.
The approximation $P_r$ is then searched in the subspace $\cU_{r}$ under the form
\begin{equation}
  P_r = \sum_{i \in \cI_{r}} \alpha_i \bigotimes_{\mu=1}^{d} Q_i^\mu.\label{eq:Pr-alpha}
\end{equation}
If the dimension $\cU_{r}$ is sufficiently small, 
$P_{r}$ is defined as the projection of $A^{-1}$ on the
subspace $\cU_r$ with respect to the inner product norm $\normt{\cdot}_{\star}$,
\begin{equation*}
  P_r = \argmin{P \in \cU_r}\norm{A^{-1}-P}_{\star},
\end{equation*}
and $P_{r}$ is given under the form \eqref{eq:Pr-alpha} with $\alpha$ solution of 
the linear system
\begin{equation*}
  \sum_{i \in \cI_{r}} \left\langle Q_{i},Q_{j} \right\rangle_{\star}\alpha_i = \left\langle A^{-1},Q_{j} \right\rangle_{\star}, \quad
  \forall j \in \cI_{r}.
\end{equation*}
When the dimension of $\cU_{r}$ is large (e.g. for large $d$ or large $r$), the computation of the projection in $\cU_{r}$ may be prohibitive. 
In this case, we replace the projection in $\cU_{r}$ by an approximate projection using Hierarchical Tucker format. More precisely, given a tree of dimensions $T$ and a set of ranks $(r_{t})_{t\in T}$, we introduce the subset of Hierarchical Tucker tensors $\cH_r^T(\cU_r)$ in $\cU_{r}$. Then, an approximation $P_{r}$ is computed by solving the minimization problem
\begin{equation*}
\min_{P \in \cH_r^T(\cU_r)} \norm{A^{-1}-P}_{\star}, 
\end{equation*}
which is equivalent to computing an approximation $P_{r}$ under the form \eqref{eq:Pr-alpha} 
with $\alpha$ in $\cH_r^T (\bigotimes_{\mu} \bR^{r_{\mu}})$ solution of 
\begin{equation*}
\min_{\alpha \in\cH_r^T (\bigotimes_{\mu} \bR^{r_{\mu}})} \norm{A^{-1}-\sum_{i \in \cI_{r}} \alpha_i \bigotimes_\mu Q_i^\mu}_{\star}.
\end{equation*}
This optimization problem is solved using an alternating minimization algorithm 
which consists in successively minimizing over the transfer tensors $\seqt{\beta_t}{t\in I(T)}$ associated with the Hierarchical Tucker representation of $\alpha$ (see Section \ref{sec:hierarchical-tensors}).  Let $\alpha = F(\{\beta^{t}\}_{t\in I(T)})$ be a parametrization of $\alpha$, where $F$ is a multilinear map, and let $F_{t}$ be the partial  application of $F$ associated with the transfer tensor $\beta^{t}$, $t\in I(T)$, such that $F_t(\beta^t) = F(\{\beta^{t}\}_{t\in I(T)})$. For a given $t\in I(T)$, the minimization on $\beta^{t}$ is equivalent to solving 
a linear system 
\begin{align*}
  N^t \beta^t = S^t,
\end{align*}
where  $N^t$ and   $S^t$ are defined such that for all tensors $\beta^{t}$ and $\delta\beta^{t}$,
\begin{align*}
\left\langle N^{t}\beta^{t} ,\delta\beta^{t} \right\rangle &=  \sum_{i \in \cI_{r}} \sum_{j \in \cI_{r}} \left\langle Q_{i},Q_{j} \right\rangle_{\star}F_{t}(\beta^{t})_i F_{t}(\delta \beta^{t})_j,
\\
\left\langle S^{t} ,\delta\beta^{t} \right\rangle&=\sum_{j \in \cI_{r}} \left\langle A^{-1},Q_{j} \right\rangle_{\star} F_{t}(\delta \beta^{t})_j   .
\end{align*}
For the practical computation of $N^{t}$ and $S^{t}$, we rely on evaluations of the inner product that exploits 
the hierarchical tensor structure (see \cite{kressner_preconditioned_2011} for technical details).
When $N^t$ is ill-conditioned, we rather define $\beta^{t}$ such that 
\begin{align*}
  \beta^t \in \arg\min_{\hat \beta^t} \norm{N^t \hat \beta^t - S^t}.
\end{align*}
The algorithm for the construction of $P_r$ is summarized in Algorithm~\ref{algo:proj_meth}.

\begin{algorithm}[H]
  \SetAlgoLined
  \KwData{$A \in \cL(\cV)$, $R \in \bN^*$, $P_{0}$}
  \KwResult{$P_R \in \cL(\cV)$}
  \For{$\mu = 1,\hdots,d$}{
      $\cU_0^\mu = 0$ \;
    }
  \For{$r=1,\hdots,R$}{
    Compute $W_r = \bigotimes_\mu W_r^\mu$ by solving $\min_{W \in
      \cC_1(\cW)} \normt{A^{-1}-P_{r-1} - W}_\star$\;
    \For{$\mu = 1,\hdots,d$}{
      $\cU_r^\mu = \cU_{r-1}^\mu + \spat{W_r^\mu}$\;
    }
    $\cU_r = \bigotimes_\mu \cU_r^\mu$\;
    Compute $P_r$ by solving $\min_{P \in \cM}\normt{A^{-1}-P}_{\star}$ with
    $\cM = \cU_r$ or $\cM = \cH_r^T(\cU_r)$\;
  }
  \Return{$P_{R}$}\;
  \caption{Projection-based algorithm for low-rank approximate inverse construction}
  \label{algo:proj_meth}
\end{algorithm}

 \section{Numerical examples}
\label{sec:applications}

In this section, we apply the proposed Algorithm \ref{algo:proj_meth} for the construction of low-rank approximations $P_{r}$ of the inverse of operators arising from discretizations of (stochastic) partial differential equations. This algorithm is denoted ALG-P. 
It is compared with a pure 
greedy rank-one algorithm for which approximations $P_{r}$ are  defined by 
 $P_{r}=P_{r-1}+W_{r}$, where the $W_{r}$ are successive rank-one corrections obtained by solving \eqref{eq:greedy}, and $P_{r}-P_{0}$ is a tensor with canonical rank $r$. This algorithm is denoted by ALG-G. For the manipulation of hierarchical Tucker format, we have used the MATLAB\textsuperscript{\textregistered} toolbox presented in \cite{kressner_htucker_2012}. 

\subsection{Poisson equation}
\label{sec:poissons-problem}

\subsubsection{Description of the problem}
\label{sec:problem-settings}

 We
consider the Poisson equation defined on the $d$-dimensional domain 
 $\Omega = \omega^d \subset \bR^d$, with $\omega = (0,1)$, whose solution $v(x)$, $x=(x_{1},\hdots,x_{d}) \in \Omega$, satisfies
\begin{align*}
  &-\sum_{\nu=1}^{d } \frac{\partial^{2} v}{\partial x_{\nu}^{2}}   = 1 \quad \text{on} \quad  \Omega \quad \text{and} \quad 
   v = 0 \quad \text{on} \quad \partial\Omega.
 \end{align*}
An approximation $u$ of  the weak solution $v\in H_0^1(\Omega) = \bigotimes_{\mu=1}^{d} H_0^1(\omega)  $ is obtained through Galerkin projection in the finite element space $\cV = \bigotimes_{\mu=1}^{d} \cV^{\mu} \subset H_0^1(\Omega)$, where the $\cV^\mu \subset H_0^1(\omega) $ are uni-dimensional linear finite element spaces. We let $\cV^\mu = \spat{\varphi_i; 1
  \le i \le n} $, where $\seqt{\varphi_i}{1 \le i
  \le n}$ is  the linear finite element basis associated with a
regular mesh of $\omega$ composed by $n+1$ elements. Therefore, 
$\cV = \{u= \sum_{i_{1}=1}^{n} \hdots \sum_{i_{d}=1}^{n} \alpha_{i_{1}\hdots i_{d}} \otimes_{\mu=1}^{d} \varphi_{i_\mu} ; 
\alpha\in \bR^{n}\otimes \hdots \otimes \bR^{n}\}$ and the Galerkin approximation 
$u\in \cV$ is defined by 
  \begin{equation}
    \sum_{\nu=1}^{d} \int_{\Omega} \frac{\partial \delta u}{\partial x_{\nu}} \frac{\partial u}{\partial x_{\nu}} \ dx =
    \int_\Omega \delta u\ dx, \quad \forall \delta u \in \cV. \label{poisson-galerkin}
\end{equation}
By identifying $u$ with the $d$-order tensor of its coefficients, also denoted
 $u\in \bR^{n}\otimes \hdots \otimes \bR^{n}$, equation \eqref{poisson-galerkin} is equivalent to the linear system 
 $Au=b$, with 
\begin{align*}
 & A = \sum_{\nu=1}^d \bigotimes_{\mu=1}^{d} \left(\delta_{\nu\mu}
    K + (1-\delta_{\nu\mu}) M\right)\ , \quad b = \bigotimes_{\mu=1}^{d} c\ , \\
&  K_{ij} = \int_{\omega}  \varphi_i'(x) \varphi_j'(x)\
  dx\ , \quad M_{ij} = \int_{\omega} \varphi_i(x)\ \varphi_j(x)\
  dx\ , \quad c_i = \int_\omega \varphi_i(x) dx\ ,
\end{align*}
and where $\delta_{\nu\mu}$ denotes the Kronecker's delta. In the following, we use $d = 20$ and $n= 100$.
\subsubsection{Numerical results}
\begin{figure}[h]
  \setlength\figureheight{0.35\textwidth}
  \setlength\figurewidth{0.5\textwidth}
  \centering
%
%
%
%
\begin{tikzpicture}

\begin{semilogyaxis}[%
view={0}{0},
width=\figurewidth,
height=\figureheight,
scale only axis,
xmin=0, xmax=10,
xlabel={r},
ymin=1e-09, ymax=0.1,
yminorticks=true,
ylabel={$\epsilon(P_{r})$},
legend style={at={(0.03,0.03)},anchor=south west,nodes=right}]
\addplot [
color=blue,
solid,
mark=+,
mark options={solid}
]
coordinates{
 (1,0.0667741451356223)(2,0.0308251334100715)(3,0.0317551783713259)(4,0.0364214003835891)(5,0.0192876399722044)(6,0.0207618088970771)(7,0.0210550661832084)(8,0.023104351065819)(9,0.022457206531031)(10,0.0229864471798287) 
};
\addlegendentry{\footnotesize ALG-G};

\addplot [
color=red,
solid,
mark=o,
mark options={solid}
]
coordinates{
 (1,0.0667741451356223)(2,0.00548979563415871)(3,0.000564660552068856)(4,6.70640036169369e-05)(5,9.6743593242561e-06)(6,2.08563722541871e-06)(7,3.52998212206799e-07)(8,7.48310758896848e-08)(9,1.43275251647949e-08)(10,2.82921586618718e-09) 
};
\addlegendentry{\footnotesize ALG-P};

\addplot [
color=blue,
dashed,
mark=+,
mark options={solid}
]
coordinates{
 (1,0.0667741451356212)(2,0.0308251334100664)(3,0.0317551783713193)(4,0.0364214003835832)(5,0.0192876399722022)(6,0.0207618088970745)(7,0.0210550661832035)(8,0.0231043510658129)(9,0.0224572065310256)(10,0.0229864471798247) 
};
\addlegendentry{\footnotesize ALG-G sym.};

\addplot [
color=red,
dashed,
mark=o,
mark options={solid}
]
coordinates{
 (1,0.0667741451356212)(2,0.00548979563415905)(3,0.000564660552067536)(4,6.70640036169251e-05)(5,9.67435074992979e-06)(6,2.08563686439918e-06)(7,3.52998037273601e-07)(8,7.48310107581902e-08)(9,1.43277545610439e-08)(10,9.46498755706829e-09) 
};
\addlegendentry{\footnotesize ALG-P sym};

\end{semilogyaxis}
\end{tikzpicture}%
  \caption{Convergence with $r$ of  the sequence of approximations $P_{r}$ computed with ALG-P (Algorithm  \ref{algo:proj_meth})  or ALG-G (pure greedy algorithm), with or without imposition of symmetry.}
  \label{fig:poisson_conv_precond}
\end{figure}
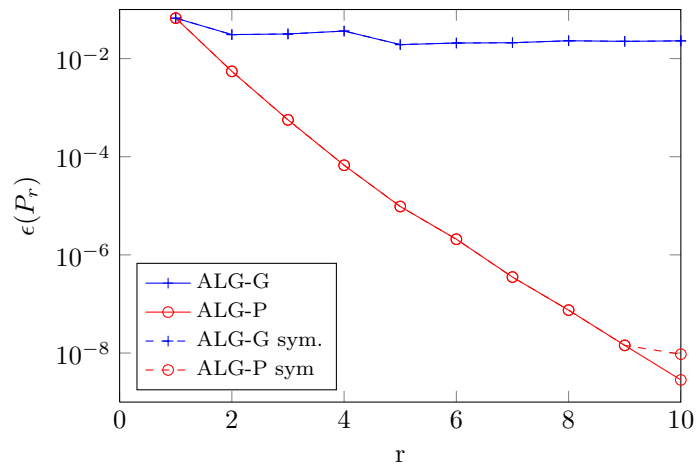

\paragraph{Convergence of the sequence of preconditioners}
The operator being symmetric
positive definite, the approximate inverse is computed with respect to
the norm $\normt{\cdot}_\star=\normt{\cdot}_A$.
For both algorithms ALG-P and ALG-G, we start from $P_{0}=0$. 
In Figure \ref{fig:poisson_conv_precond} we observe the convergence of
the approximate inverse $P_r$, by computing the relative error
\begin{equation}\label{eq:estim_err}
\epsilon(P_{r}) = \frac{\norm{I-P_r A}}{\norm{I}} =
\frac{\norm{A^{-1}-P_r}_{AA^T}}{\norm{A^{-1}}_{AA^T}}.
\end{equation}
The convergence with $r$ of $\epsilon(P_{r})$ is plotted for both algorithms, with or without imposition of symmetry in the rank-one corrections.  
For algorithms without imposition of symmetry, we  first observe the convergence of $P_{r}$ towards 
 $A^{-1}$. ALG-P provides a sequence of approximations $P_{r}$ which converges very fast compared to the greedy construction. For the same rank $r=10$, ALG-P and ALG-G yield errors of $3.10^{-9}$ and $3.10^{-2}$ respectively.
Algorithms with imposition of symmetry present almost the same convergence, 
except for a rank greater than 9 corresponding to a very  
small error less than $2.10^{-8}$. 
Let us emphasize  that for computing $P_{r}$, both algorithms require the computation of the  same number $r$ of rank-one corrections. For large matrices, these rank-one corrections constitute the most costly step of this algorithm and therefore, for computing $P_{r}$, the two algorithms
require almost the same computation times. 

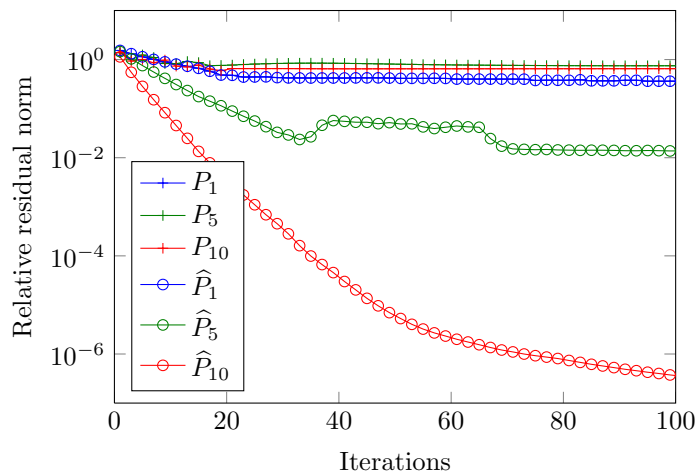
\begin{figure}[h]
  \centering
  \setlength\figureheight{0.35\textwidth}
  \setlength\figurewidth{0.5\textwidth}
%
%
%
%
\begin{tikzpicture}

\begin{semilogyaxis}[%
width=\figurewidth,
height=\figureheight,
scale only axis,
xmin=0,
xmax=100,
xlabel={Iterations},
ymin=1e-07,
ymax=10,
yminorticks=true,
ylabel={Relative residual norm},
legend style={at={(0.03,0.03)},anchor=south west,nodes=right}
]
\addplot [
color=blue,
solid,
mark=+,
mark options={solid}
]
table[row sep=crcr]{
1 1.53677320945239\\
3 1.32004286327684\\
5 1.14072373229854\\
7 1.00823212334204\\
9 0.89802386632799\\
11 0.808354962355744\\
13 0.732050164936217\\
15 0.663618529602625\\
17 0.578483597291053\\
19 0.483947380396591\\
21 0.478276899164633\\
23 0.440637292448813\\
25 0.442561682267017\\
27 0.445148138478953\\
29 0.429372568535783\\
31 0.423298327792494\\
33 0.422144378894488\\
35 0.421176679699965\\
37 0.421248035140315\\
39 0.421504393128299\\
41 0.420758771217229\\
43 0.423553626399056\\
45 0.423931159859476\\
47 0.424079762302746\\
49 0.420344874953542\\
51 0.413876079419341\\
53 0.414207542266996\\
55 0.41384756155984\\
57 0.415564076937687\\
59 0.413127783220415\\
61 0.400343610790663\\
63 0.40042877642742\\
65 0.400373078320771\\
67 0.40035694409873\\
69 0.400150187480588\\
71 0.40012937782353\\
73 0.398699464990996\\
75 0.375626742983772\\
77 0.37991674495901\\
79 0.378794812219772\\
81 0.378813779408778\\
83 0.387773375577188\\
85 0.36952931804963\\
87 0.3692253586489\\
89 0.369079443279251\\
91 0.378621761238379\\
93 0.378753612815989\\
95 0.362675813384552\\
97 0.362646022485589\\
99 0.362569496059219\\
101 0.363185469570729\\
};
\addlegendentry{$P_1$};

\addplot [
color=green!50!black,
solid,
mark=+,
mark options={solid}
]
table[row sep=crcr]{
1 1.37541505252121\\
3 0.992672081755847\\
5 1.24380880864479\\
7 1.10230683146114\\
9 0.831473223467727\\
11 0.826488147866066\\
13 0.912371208487396\\
15 0.868780457612802\\
17 0.752195774331185\\
19 0.759608771869494\\
21 0.773726915411443\\
23 0.793505176471371\\
25 0.810300937395698\\
27 0.823156522198397\\
29 0.833158546800424\\
31 0.840283649318114\\
33 0.844536765848065\\
35 0.846096987475455\\
37 0.845255286515211\\
39 0.842377067934385\\
41 0.837868724028382\\
43 0.832144646713097\\
45 0.825600124274999\\
47 0.818595222440344\\
49 0.811447028241415\\
51 0.804422772474359\\
53 0.79773194547326\\
55 0.791546878753652\\
57 0.786062907573693\\
59 0.781245096276381\\
61 0.777029772825865\\
63 0.77334810848495\\
65 0.770089750191462\\
67 0.767146147333443\\
69 0.764441428812135\\
71 0.761934636868076\\
73 0.75961159508824\\
75 0.757475883945048\\
77 0.755541004072533\\
79 0.753823931757174\\
81 0.752340177356266\\
83 0.751100370726469\\
85 0.750108169041973\\
87 0.749359203462296\\
89 0.748840875382383\\
91 0.748532909651211\\
93 0.748408576995554\\
95 0.748436440228741\\
97 0.748582419026139\\
99 0.748811940628737\\
101 0.7490919561324\\
};
\addlegendentry{$P_5$};

\addplot [
color=red,
solid,
mark=+,
mark options={solid}
]
table[row sep=crcr]{
1 1.38502686567431\\
3 0.930978086163214\\
5 1.00714033892595\\
7 0.859945099995247\\
9 0.978043434353648\\
11 0.736297160364572\\
13 0.708302715181589\\
15 0.774153788020624\\
17 0.641649060502777\\
19 0.59561423943131\\
21 0.655142939074308\\
23 0.654760257762895\\
25 0.65420332654516\\
27 0.654732046222349\\
29 0.655332393236842\\
31 0.654940904395761\\
33 0.6535817260425\\
35 0.651632449781648\\
37 0.649858985908741\\
39 0.648768369707993\\
41 0.648399251397025\\
43 0.648522658301373\\
45 0.648824114028951\\
47 0.649099108516847\\
49 0.649278364014851\\
51 0.649363183116141\\
53 0.64938374457571\\
55 0.649379736615335\\
57 0.649385371737821\\
59 0.649421061770356\\
61 0.649494739283425\\
63 0.649607899008599\\
65 0.649759876045643\\
67 0.649947753286506\\
69 0.650163064775273\\
71 0.650390488571045\\
73 0.650612331194754\\
75 0.650815177441248\\
77 0.650993320364292\\
79 0.65114860032492\\
81 0.651288382759945\\
83 0.651422499626631\\
85 0.651559842384474\\
87 0.651705496670641\\
89 0.651859067865743\\
91 0.652014410322082\\
93 0.652160750471921\\
95 0.652285137830958\\
97 0.652375867139805\\
99 0.652425973590982\\
101 0.652435671327349\\
};
\addlegendentry{$P_{10}$};

\addplot [
color=blue,
solid,
mark=o,
mark options={solid}
]
table[row sep=crcr]{
1 1.53677320945239\\
3 1.32004286327684\\
5 1.14072373229854\\
7 1.00823212334204\\
9 0.89802386632799\\
11 0.808354962355744\\
13 0.732050164936217\\
15 0.663618529602625\\
17 0.578483597291053\\
19 0.483947380396591\\
21 0.478276899164633\\
23 0.440637292448813\\
25 0.442561682267017\\
27 0.445148138478953\\
29 0.429372568535783\\
31 0.423298327792494\\
33 0.422144378894488\\
35 0.421176679699965\\
37 0.421248035140315\\
39 0.421504393128299\\
41 0.420758771217229\\
43 0.423553626399056\\
45 0.423931159859476\\
47 0.424079762302746\\
49 0.420344874953542\\
51 0.413876079419341\\
53 0.414207542266996\\
55 0.41384756155984\\
57 0.415564076937687\\
59 0.413127783220415\\
61 0.400343610790663\\
63 0.40042877642742\\
65 0.400373078320771\\
67 0.40035694409873\\
69 0.400150187480588\\
71 0.40012937782353\\
73 0.398699464990996\\
75 0.375626742983772\\
77 0.37991674495901\\
79 0.378794812219772\\
81 0.378813779408778\\
83 0.387773375577188\\
85 0.36952931804963\\
87 0.3692253586489\\
89 0.369079443279251\\
91 0.378621761238379\\
93 0.378753612815989\\
95 0.362675813384552\\
97 0.362646022485589\\
99 0.362569496059219\\
101 0.363185469570729\\
};
\addlegendentry{$\widehat P_{1}$};

\addplot [
color=green!50!black,
solid,
mark=o,
mark options={solid}
]
table[row sep=crcr]{
1 1.45417492713106\\
3 1.03538659217619\\
5 0.755139606988628\\
7 0.557066445897503\\
9 0.414346821163279\\
11 0.3102660467064\\
13 0.233555135605115\\
15 0.176396690224945\\
17 0.142096212581118\\
19 0.114501865565083\\
21 0.0891571110636637\\
23 0.0700042266688918\\
25 0.0558538489160291\\
27 0.0429559760666546\\
29 0.0328585408301081\\
31 0.0291486682889387\\
33 0.0235214199222556\\
35 0.026710680839776\\
37 0.0459136753111817\\
39 0.0572196062048067\\
41 0.0551404082379756\\
43 0.0529280235636172\\
45 0.0528676215028746\\
47 0.0493869252801382\\
49 0.0519175358572195\\
51 0.0487203680648099\\
53 0.0491858916120133\\
55 0.0425442311710622\\
57 0.0393330615367106\\
59 0.0417993210721111\\
61 0.0443443768683927\\
63 0.0425260384379452\\
65 0.0417781068669025\\
67 0.0245425071465739\\
69 0.0172245504807731\\
71 0.0151135406045894\\
73 0.0147943027949734\\
75 0.0145675946921165\\
77 0.0147415576728849\\
79 0.0144029175674753\\
81 0.0142730373664109\\
83 0.0141451556337668\\
85 0.0142219773743176\\
87 0.0142082251790218\\
89 0.0140908127229912\\
91 0.0138990101173277\\
93 0.0138675040590285\\
95 0.0140300964985843\\
97 0.0139262897422533\\
99 0.0137522086365811\\
101 0.0137782520886972\\
};
\addlegendentry{$\widehat P_{5}$};

\addplot [
color=red,
solid,
mark=o,
mark options={solid}
]
table[row sep=crcr]{
1 1.13646316121122\\
3 0.551084913445603\\
5 0.284192998515266\\
7 0.152072179770469\\
9 0.0827647843766511\\
11 0.0453463602476459\\
13 0.0248176557702889\\
15 0.0134687017123422\\
17 0.00785282402069343\\
19 0.00477335884165926\\
21 0.002674884378089\\
23 0.00174431476201406\\
25 0.00110057009929271\\
27 0.000689335048131677\\
29 0.000452480640892707\\
31 0.000281402000352129\\
33 0.000161503257357482\\
35 9.92020307458686e-05\\
37 6.65228125954341e-05\\
39 4.56032862426295e-05\\
41 3.00203039813245e-05\\
43 2.02677346943952e-05\\
45 1.36284102442708e-05\\
47 9.59725284294854e-06\\
49 6.96835459501313e-06\\
51 5.2824734376376e-06\\
53 4.02386643899015e-06\\
55 3.22376509144237e-06\\
57 2.69190699549363e-06\\
59 2.32314103234745e-06\\
61 1.99821771700664e-06\\
63 1.7553476993138e-06\\
65 1.5455076627241e-06\\
67 1.35679335333818e-06\\
69 1.2061164276666e-06\\
71 1.10238030738646e-06\\
73 9.98980305971454e-07\\
75 9.24060625787806e-07\\
77 8.65384064301849e-07\\
79 8.08974917753999e-07\\
81 7.38201330995064e-07\\
83 6.74868318074066e-07\\
85 6.16920120374804e-07\\
87 5.70910618946884e-07\\
89 5.28564329193504e-07\\
91 4.8805546968289e-07\\
93 4.56894931842127e-07\\
95 4.2554204065887e-07\\
97 4.00721131333299e-07\\
99 3.74496143025722e-07\\
101 3.54624481196449e-07\\
};
\addlegendentry{$\widehat P_{10}$};

\end{semilogyaxis}
\end{tikzpicture}%
  \caption{Convergence of the Preconditioned Conjugate Gradient using approximate iterates in $\cH_{15}^T(\cV)$, and using preconditioners $\widehat P_{r}$ (resp. $P_{r}$) constructed using ALG-P   (resp. 
  ALG-G) with different  ranks $r \in \{1,5,10\}$.}
  \label{fig:poisson_conv_sol}
\end{figure}

\paragraph{Preconditioned iterative solver}
The operator $A$ being symmetric positive definite, we solve the
linear system $Au=b$ using a Preconditioned Conjugate Gradient (PCG)
with low-rank tensor compressions of the iterates in hierarchical
Tucker format (see the algorithm in \cite{kressner_low-rank_2011}).
Here, we use approximations of the iterates in the hierarchical Tucker
subset $\cH_{15}^T(\cV)$ associated with a balanced tree $T$ (same
rank $15$ at each node of the tree).  We analyze the convergence of
the PCG using symmetric preconditioners $P_{r}$ constructed with ALG-P
or ALG-G.  On Figure \ref{fig:poisson_conv_sol}, we observe that the
convergence rate of the PCG strongly depends on the quality of the
preconditioner. We first note that when using a preconditioner $P_{r}$
constructed with a pure greedy algorithm, the convergence of the PCG
is not improved when increasing the rank $r$ of the preconditioner.
However, when using ALG-P, we can see that the convergence rate
rapidly increases with the rank $r$.  Moreover, we observe
that the relative residual norm stagnates at a certain precision. This precision depends on two factors: the low-rank subset chosen for the approximation of iterates (here fixed at $\cH_{15}^T(\cV)$) and the quality of the preconditioner. Figure \ref{fig:poisson_conv_sol} illustrates the strong influence of the preconditioner on the resulting precision, and the superiority of the proposed algorithm over the pure greedy construction. 
In particular, we observe a difference of 2 (resp. 7) orders of magnitude between rank-$5$ (resp. rank-$10$) preconditioners constructed by ALG-G and ALG-P. 

\subsection{Stochastic elliptic problem}
\label{sec:elliptic-equation}

\subsubsection{Description of the problem}
\label{sec:discr-param-probl}

We consider the stochastic partial differential equation defined on a 2-dimensional domain $\Omega=(0,1)^{2}$,
\begin{align*}
-\kappa \Delta v + \eta  v = f\quad \text{on} \quad \Omega, \quad v=0\quad \text{on} \quad \partial  \Omega,
\end{align*}
where $\kappa$ and $\eta$ are independent random variables with
uniform law over the intervals $(1,10)$ and $(200,1000)$
respectively. $f$ is such that $f(x) = 1$ if $x \in [0.6,0.8] \times
[0.6,0.8]$ and $f(x) = 0$ otherwise. $\kappa$ (resp. $\eta$) is
expressed as a linear function $\kappa(\xi_{1})$
(resp. $\eta(\xi_{2})$ of a random variable $\xi_{1}$
(resp. $\xi_{2})$ with uniform law on the interval $\Xi^{1}=(-1,1)$
(resp. $\Xi^{2}$), and the solution is expressed under the form
$v(x,\xi_{1},\xi_{2})$, $x\in \Omega$, with $v: \Omega\times \Xi^{1}
\times \Xi^{2} \rightarrow \bR$.  A Galerkin approximation $u$ of the
weak solution $v \in H^1_0(\Omega) \otimes L^2(\Xi^{1}) \otimes
L^2(\Xi^{2})$ is obtained though Galerkin projection in a finite
dimensional space $\cV = \cV^{1} \otimes \cV^{2} \otimes \cV^{3}$,
where $\cV^1 = \spat{\varphi_i;\ 1 \le i \le n} \subset H_0^1(\Omega)$
is a finite element space associated with a regular cartesian mesh of
$\Omega$, and where $\cV^{2}\subset L^2(\Xi^{1})$ and $\cV^{3}\subset
L^2(\Xi^{2})$ are polynomial spaces of degree $p-1$. We denote by
$\{\varphi_{i},1\le i\le n\}$ the finite element basis of $\cV^{1}$.
For $\cV^{2}$ and $\cV^{3}$, we introduce a basis $\{\psi_{i};1\le i
\le p\}$ where $\psi_{i+1}$ denotes the Legendre polynomial of degree
$i$.  Therefore, $\cV = \{u= \sum_{i=1}^{n} \sum_{j=1}^{p}
\sum_{k=1}^{p} \alpha_{ijk} \varphi_{i}\otimes \psi_{j }\otimes
\psi_{k} ; \alpha \in \bR^{n\times p\times p}\}$ and the Galerkin approximation $u\in \cV$ is defined by
\begin{equation} \label{eq:elliptic-pb}
  \int_{\Omega \times \Xi^{1} \times \Xi^{2}} (\kappa(y_{1}) \nabla \delta v
  \cdot  \nabla v + \eta(y_{2}) \delta v~v)~dxdy_{1}dy_{2} =
  \int_{\Omega \times \Xi^{1} \times \Xi^{2}} \delta v~f(x) ~dxdy_{1}dy_{2}
\end{equation}
for all $\delta v \in \cV$. By identifying $u$ with the $3$-order tensor of its coefficients, also denoted
 $u\in \bR^{n}\otimes \bR^{p} \otimes \bR^{p} $, equation \eqref{eq:elliptic-pb} is equivalent to a linear system 
 $Au=b$, where $A$ is a rank-2 operator and $b$ is a rank-one tensor. 
In the following numerical experiments, we take 
$n= 20^2$ and $p = 10$.

\subsubsection{Numerical results}
\label{sec:numerical-results}

\paragraph{Convergence of the sequence of preconditioners}

\begin{figure}[h]
  \setlength\figureheight{0.35\textwidth}
  \setlength\figurewidth{0.5\textwidth}
  \centering
%
%
%
%
\begin{tikzpicture}

\begin{semilogyaxis}[%
view={0}{0},
width=\figurewidth,
height=\figureheight,
scale only axis,
xmin=1, xmax=10,
xlabel={$r$},
ymin=1e-08, ymax=1,
yminorticks=true,
ylabel={$\epsilon(P_{r})$},
legend style={at={(0.03,0.03)},anchor=south west,nodes=right}]
\addplot [
color=blue,
solid,
mark=+,
mark options={solid}
]
coordinates{
 (1,0.141454782530368)(2,0.0884997770703733)(3,0.0715455620394111)(4,0.0424634168330422)(5,0.0392779813857653)(6,0.0188698876933793)(7,0.0144582778939564)(8,0.0100577985425702)(9,0.00770517432767674)(10,0.00765771127322665) 
};
\addlegendentry{ALG-G, $\gamma=100\%$};

\addplot [
color=red,
solid,
mark=o,
mark options={solid}
]
coordinates{
 (1,0.140076700333146)(2,0.102380449176266)(3,0.0789211308779267)(4,0.063194076602968)(5,0.0439352319779074)(6,0.0424899027162276)(7,0.0367453811189302)(8,0.0365134590086583)(9,0.0333820741669001)(10,0.0318097945537572) 
};
\addlegendentry{ALG-G, $\gamma=10\%$};

\addplot [
color=green!50!black,
solid,
mark=square,
mark options={solid}
]
coordinates{
 (1,0.140576760631264)(2,0.0873785366295442)(3,0.0546834655859887)(4,0.0480702952905874)(5,0.0238265156886309)(6,0.0207240109411102)(7,0.0189824614871396)(8,0.0182087177074357)(9,0.0167997798820887)(10,0.0123962933919746) 
};
\addlegendentry{ALG-G, $\gamma = 30\%$};

\addplot [
color=blue,
dashed,
mark=+,
mark options={solid}
]
coordinates{
 (1,0.141454782530368)(2,0.0425791221955531)(3,0.0136887787368225)(4,0.00370179233427925)(5,0.000634801621685698)(6,0.000120213642725391)(7,8.52407522199584e-06)(8,1.35429492377324e-06)(9,1.75346998538299e-07)(10,5.29900029886889e-08) 
};
\addlegendentry{ALG-P, $\gamma=100\%$};

\addplot [
color=red,
dashed,
mark=o,
mark options={solid}
]
coordinates{
 (1,0.140076700333146)(2,0.0632210883561922)(3,0.0546971892585285)(4,0.0390995988407028)(5,0.0337375571293828)(6,0.0299253404363049)(7,0.0264377409203635)(8,0.0241441384567534)(9,0.0221143415511827)(10,0.0209057571778937) 
};
\addlegendentry{ALG-P, $\gamma=10\%$};

\addplot [
color=green!50!black,
dashed,
mark=square,
mark options={solid}
]
coordinates{
 (1,0.140576760631264)(2,0.0446917590597903)(3,0.0207329330788663)(4,0.0141461991033271)(5,0.0107384360888911)(6,0.0095334399226097)(7,0.00835768590497494)(8,0.00730745707616945)(9,0.00664296707553744)(10,0.00616146591478597) 
};
\addlegendentry{ALG-P, $\gamma = 30\%$};

\end{semilogyaxis}
\end{tikzpicture}%
  \caption{Convergence with $r$ of  the sequence of approximations $P_{r}$ computed with ALG-P (Algorithm  \ref{algo:proj_meth})  or ALG-G (pure greedy rank-one algorithm). Sparsity along dimension $1$ is imposed with a maximum percentage $\gamma$ of nonzero components, for $\gamma = 10,30,100 \%$. }
  \label{fig:largeleaf_conv_precond}
\end{figure}
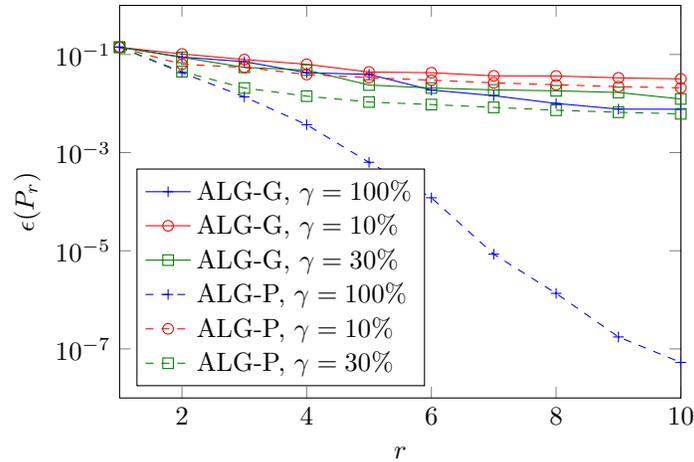

We construct sequences of low-rank preconditioner $P_{r}$ using either
ALG-P (Algorithm \ref{algo:proj_meth}) or ALG-G (pure greedy
algorithm), with a norm $\normt{\cdot}_\star=\normt{\cdot}_{A}$.  We
impose sparsity only along dimension 1. More precisely, during the
computation of a rank-one correction $W_{r}=W_{r}^{1}\otimes
W_{r}^{2}\otimes W_{r}^{3}$, we search for a sparse approximation
$W_{r}^{1}$ using the adaptive algorithm presented in Section
\ref{sec:imposing-sparsity}, with a number of nonzero components
limited to a certain percentage, denoted $\gamma$, of the total number
of components $n^{2}$.  The convergence with $r$ of the different
preconditioners is illustrated in Figure
\ref{fig:largeleaf_conv_precond}, where the error estimator
$\epsilon(P_{r})$ is defined by equation \eqref{eq:estim_err}.  We
observe that all algorithms seem to converge toward the inverse of the
operator $A$. As we could have expected, the higher $\gamma$, the
better the approximation is. Compared to ALG-G, ALG-P significantly
improves the quality of the approximation for a low additional
cost. For $\gamma = 30\%$, ALG-P provides for $r=4$ a better
approximation than the approximation provided by ALG-G  for $r=10$.

\paragraph{Preconditioned iterative solver}
\label{sec:impr-due-proj}

Although the operator $A$ is symmetric positive definite, the
resulting preconditioner with imposed sparsity is non symmetric. For
the solution of the linear system, we therefore use a preconditioned
GMRES (without restart) with low-rank approximations of the iterates
(see the algorithm in \cite{Ballani2013}). Approximations of the
iterates are here computed in the subset $\cT_{m}(\cV)$ of
rank-$(m,m,m)$ Tucker tensors using the Higher Order Orthogonal
Iterations (HOOI) algorithm \cite{DeLathauwer2000}.

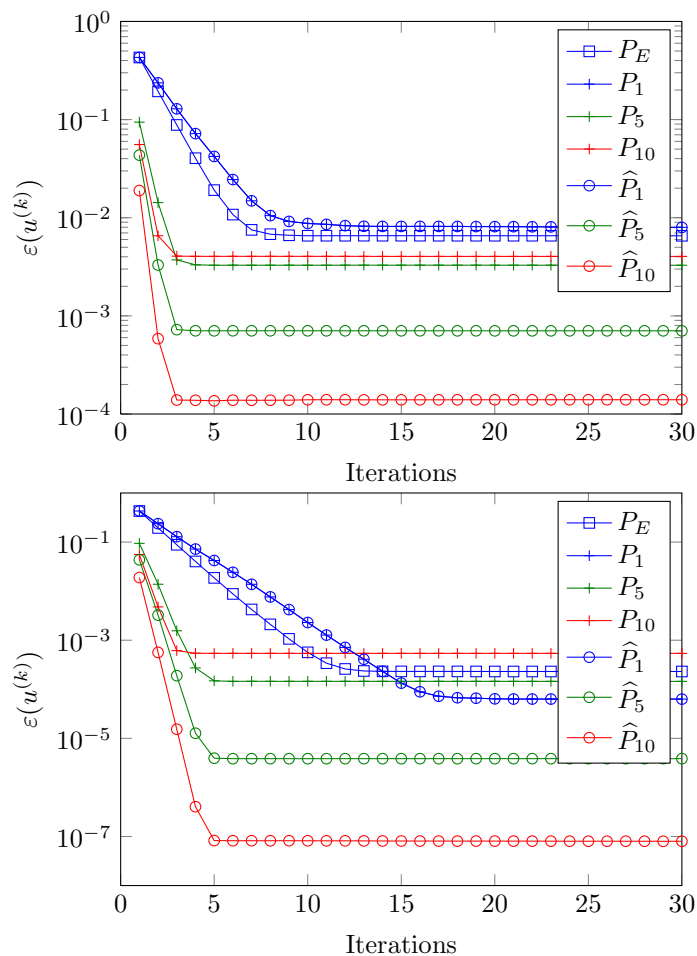
\begin{figure}[h]
  \centering
  \setlength\figureheight{0.35\textwidth}
  \setlength\figurewidth{0.5\textwidth}
%
%
%
\begin{tikzpicture}

\begin{semilogyaxis}[%
width=\figurewidth,
height=\figureheight,
scale only axis,
xmin=0,
xmax=30,
xlabel={Iterations},
ymin=0.0001,
ymax=1,
yminorticks=true,
ylabel={$\varepsilon(u^{(k)})$},
legend style={draw=black,fill=white,legend cell align=left}
]
\addplot [
color=blue,
solid,
mark=square,
mark options={solid}
]
table[row sep=crcr]{
1 0.429158943473958\\
2 0.193924386555288\\
3 0.0881700553015909\\
4 0.0404095205676266\\
5 0.0191178928523393\\
6 0.0107669391559102\\
7 0.00752185987776282\\
8 0.0067955268152443\\
9 0.00662598662639562\\
10 0.00655272277776467\\
11 0.00655684465374879\\
12 0.00655630300154933\\
13 0.0065619491578679\\
14 0.00655953836048073\\
15 0.00656075779507579\\
16 0.00656073565122061\\
17 0.00655877837983758\\
18 0.00655693407536949\\
19 0.00655558513004954\\
20 0.00655546196484764\\
21 0.00655517403767528\\
22 0.00655517814663464\\
23 0.00655510429463314\\
24 0.00655517506253689\\
25 0.00655535779451559\\
26 0.00655429425936867\\
27 0.00655329365098749\\
28 0.00655058732142624\\
29 0.00654840847902228\\
30 0.00654643683192794\\
};
\addlegendentry{$P_E$};

\addplot [
color=blue,
solid,
mark=+,
mark options={solid}
]
table[row sep=crcr]{
1 0.42992843706819\\
2 0.236781091471918\\
3 0.128670509484902\\
4 0.0719687966549017\\
5 0.042131004380971\\
6 0.024491028603022\\
7 0.0148732302433656\\
8 0.0105147016028054\\
9 0.00917398726832181\\
10 0.00874480903476679\\
11 0.00856069625972374\\
12 0.00830798554908\\
13 0.00818421994861361\\
14 0.00816184273487837\\
15 0.00815864384970394\\
16 0.00816401187895246\\
17 0.00815976880554773\\
18 0.00815201773357234\\
19 0.00813037119186941\\
20 0.00812455141007925\\
21 0.00811170181003298\\
22 0.00808894256520563\\
23 0.00806997702683564\\
24 0.00802090934389001\\
25 0.00800975697202343\\
26 0.00800644304952495\\
27 0.00798158809492071\\
28 0.00798310713606205\\
29 0.00798150771356039\\
30 0.00798073966975526\\
};
\addlegendentry{$P_1$};

\addplot [
color=green!50!black,
solid,
mark=+,
mark options={solid}
]
table[row sep=crcr]{
1 0.0941624671358326\\
2 0.0142771640885178\\
3 0.00372001721389246\\
4 0.0033121753754568\\
5 0.00327970067801952\\
6 0.00328475303478438\\
7 0.0032832797499173\\
8 0.00328318106738159\\
9 0.00328324043748591\\
10 0.00328282637300127\\
11 0.0032832537571923\\
12 0.00328497582476427\\
13 0.00328492085746962\\
14 0.00328502664573372\\
15 0.00328577344950372\\
16 0.00328626756786768\\
17 0.00328536912228475\\
18 0.00328505457767686\\
19 0.00328496522913585\\
20 0.0032845363238357\\
21 0.0032841800522554\\
22 0.00328395465268379\\
23 0.00328418266068468\\
24 0.00328418656887807\\
25 0.00328414203166666\\
26 0.00328415471081687\\
27 0.00328402116457207\\
28 0.00328401615870062\\
29 0.00328369306635182\\
30 0.00328104041220835\\
};
\addlegendentry{$P_5$};

\addplot [
color=red,
solid,
mark=+,
mark options={solid}
]
table[row sep=crcr]{
1 0.0556955195383315\\
2 0.00655061148149525\\
3 0.00406226286321615\\
4 0.00404406773584782\\
5 0.00404108967391761\\
6 0.00404385409669075\\
7 0.0040428691167917\\
8 0.00404302743469598\\
9 0.0040425770107057\\
10 0.00404251396266575\\
11 0.0040424531555702\\
12 0.0040424351176822\\
13 0.00404112101761145\\
14 0.00403715747167926\\
15 0.00403715709119155\\
16 0.00403708590799767\\
17 0.00403709743753427\\
18 0.00403710736864559\\
19 0.00403706941535468\\
20 0.00403666313764792\\
21 0.0040365420766409\\
22 0.00403667409445501\\
23 0.00403673033173021\\
24 0.00403666486965566\\
25 0.00403675289578223\\
26 0.00403638389275186\\
27 0.00403642758422167\\
28 0.00403620929057795\\
29 0.00403612970978229\\
30 0.00403607185446815\\
};
\addlegendentry{$P_{10}$};

\addplot [
color=blue,
solid,
mark=o,
mark options={solid}
]
table[row sep=crcr]{
1 0.42992843706819\\
2 0.236781091471918\\
3 0.12867050948491\\
4 0.0719687966549024\\
5 0.0421310043809724\\
6 0.0244910337041273\\
7 0.0148731898108975\\
8 0.0105114575304236\\
9 0.00916870679004603\\
10 0.00873982097691359\\
11 0.00855671875534857\\
12 0.00830494628678802\\
13 0.00818142353525232\\
14 0.00815868809939376\\
15 0.00815536381089566\\
16 0.00816075985496901\\
17 0.00815645442689757\\
18 0.00814843116304996\\
19 0.00812652318697484\\
20 0.00812071355302719\\
21 0.00810856022891975\\
22 0.0080854102847881\\
23 0.00806609530727728\\
24 0.00801697408690945\\
25 0.00800587018238534\\
26 0.00800248837898055\\
27 0.00797738569143126\\
28 0.00797891651693714\\
29 0.00797729559171634\\
30 0.00797659504210993\\
};
\addlegendentry{$\widehat P_1$};

\addplot [
color=green!50!black,
solid,
mark=o,
mark options={solid}
]
table[row sep=crcr]{
1 0.043410621253962\\
2 0.00330019657469539\\
3 0.000723833573895744\\
4 0.000705244367550545\\
5 0.000704221134906199\\
6 0.000705242839717393\\
7 0.000705131810274106\\
8 0.000705152618076653\\
9 0.000705819758970892\\
10 0.000704959465422618\\
11 0.000704805427053817\\
12 0.000704802208786707\\
13 0.000704775513261455\\
14 0.000704752192235011\\
15 0.000704666615769654\\
16 0.000704607600163609\\
17 0.00070458388693336\\
18 0.000704302329437687\\
19 0.000704946131834381\\
20 0.000704093696626489\\
21 0.000704099858543493\\
22 0.000704191466609824\\
23 0.000704190110123061\\
24 0.000704192003631242\\
25 0.000704204709636677\\
26 0.000704238572511431\\
27 0.000704242349805227\\
28 0.000704244405791343\\
29 0.000704264397982547\\
30 0.000704218987291762\\
};
\addlegendentry{$\widehat P_5$};

\addplot [
color=red,
solid,
mark=o,
mark options={solid}
]
table[row sep=crcr]{
1 0.0189227860236283\\
2 0.000585883642294805\\
3 0.000138862933447314\\
4 0.000137995477959175\\
5 0.000136289279131353\\
6 0.000138383346891256\\
7 0.000138045544062411\\
8 0.000138169504448673\\
9 0.000138368530773669\\
10 0.000139635901674672\\
11 0.000139799285053435\\
12 0.000139812492137298\\
13 0.000139468253682434\\
14 0.000139499579599182\\
15 0.000139505933327069\\
16 0.000139526733404065\\
17 0.00013952216339755\\
18 0.000139615836784086\\
19 0.000139623051548641\\
20 0.000139669956556833\\
21 0.00013967006022306\\
22 0.000139673551438144\\
23 0.000139692198289339\\
24 0.000139753839882805\\
25 0.000139781236810756\\
26 0.000139782338004014\\
27 0.000139755988357834\\
28 0.000139746068819785\\
29 0.000139771515158977\\
30 0.000139779579458205\\
};
\addlegendentry{$\widehat P_{10}$};

\end{semilogyaxis}
\end{tikzpicture}%
%
%
%
\begin{tikzpicture}

\begin{semilogyaxis}[%
width=\figurewidth,
height=\figureheight,
scale only axis,
xmin=0,
xmax=30,
xlabel={Iterations},
ymin=1e-08,
ymax=1,
yminorticks=true,
ylabel={$\varepsilon(u^{(k)})$},
legend style={draw=black,fill=white,legend cell align=left}
]
\addplot [
color=blue,
solid,
mark=square,
mark options={solid}
]
table[row sep=crcr]{
1 0.429158943473958\\
2 0.193924386555288\\
3 0.0881738905048614\\
4 0.0402423050333519\\
5 0.0186255542858626\\
6 0.00876841412129354\\
7 0.00423273501794931\\
8 0.00210268081903425\\
9 0.00106472684159046\\
10 0.000565803043768087\\
11 0.000340255257410098\\
12 0.00025934634335938\\
13 0.000236029173595588\\
14 0.000232000944803067\\
15 0.000231383921007561\\
16 0.000231254029183732\\
17 0.000231398978916289\\
18 0.000231324633049877\\
19 0.000230752141603198\\
20 0.0002305351068571\\
21 0.00023033912553766\\
22 0.000230343964222403\\
23 0.00023039444532644\\
24 0.000230370464232207\\
25 0.000230569271976058\\
26 0.000230652621298116\\
27 0.000230614400477315\\
28 0.000230495853228814\\
29 0.000230468086111631\\
30 0.000230545288526876\\
};
\addlegendentry{$P_E$};

\addplot [
color=blue,
solid,
mark=+,
mark options={solid}
]
table[row sep=crcr]{
1 0.42992843706819\\
2 0.236781091471918\\
3 0.128670455030046\\
4 0.0719661177495313\\
5 0.0421006325937604\\
6 0.024192110687207\\
7 0.0137819430932098\\
8 0.00762363762642118\\
9 0.00420728748891599\\
10 0.00229321995674468\\
11 0.001272051679938\\
12 0.00071351298972875\\
13 0.000409319228710464\\
14 0.000232312203940179\\
15 0.000133454836138397\\
16 8.89023979325352e-05\\
17 7.20770647965562e-05\\
18 6.7242053843388e-05\\
19 6.57633014154542e-05\\
20 6.39439491508067e-05\\
21 6.32827613598549e-05\\
22 6.32415308952603e-05\\
23 6.32278162861285e-05\\
24 6.32321021690902e-05\\
25 6.32419786339356e-05\\
26 6.3158627405285e-05\\
27 6.31831608078594e-05\\
28 6.31872509978776e-05\\
29 6.32690346766748e-05\\
30 6.30893033112665e-05\\
};
\addlegendentry{$P_1$};

\addplot [
color=green!50!black,
solid,
mark=+,
mark options={solid}
]
table[row sep=crcr]{
1 0.0941624671358326\\
2 0.0138278613212396\\
3 0.00155807998304732\\
4 0.000272730449180767\\
5 0.000148748150349933\\
6 0.000145062786862092\\
7 0.000145540913316611\\
8 0.000145215727061201\\
9 0.000145240195856081\\
10 0.000145242660665751\\
11 0.000145249901972283\\
12 0.000145246889494005\\
13 0.000145272567795389\\
14 0.000145272019669065\\
15 0.000145278173303531\\
16 0.000145267595021659\\
17 0.000145264267351227\\
18 0.000145259008968226\\
19 0.000145186322155119\\
20 0.000145275053869801\\
21 0.000145088231526027\\
22 0.000145085372168219\\
23 0.000145023396780709\\
24 0.00014495886842974\\
25 0.000144965394401241\\
26 0.000144976686980868\\
27 0.000144973707293241\\
28 0.000145003199975725\\
29 0.000144964992999729\\
30 0.000144985634010638\\
};
\addlegendentry{$P_5$};

\addplot [
color=red,
solid,
mark=+,
mark options={solid}
]
table[row sep=crcr]{
1 0.0557085029412057\\
2 0.0047751134433319\\
3 0.000613655575112465\\
4 0.000542365823957377\\
5 0.000539602475123364\\
6 0.000540157626737206\\
7 0.000540049381627121\\
8 0.000539982413607281\\
9 0.000540002577888536\\
10 0.000540031434341188\\
11 0.000540001859976022\\
12 0.000539881814239134\\
13 0.000539814467674161\\
14 0.000539817429720723\\
15 0.000539853224922269\\
16 0.000539855057421629\\
17 0.00053975604150994\\
18 0.00053967937309339\\
19 0.000539631816897017\\
20 0.000539576911410428\\
21 0.000539590059447635\\
22 0.000539567730336878\\
23 0.000539565822388189\\
24 0.000539548139838826\\
25 0.000539410185558509\\
26 0.00053942265410943\\
27 0.000539421153636899\\
28 0.000539411624625555\\
29 0.000539433753655577\\
30 0.000539422672265643\\
};
\addlegendentry{$P_{10}$};

\addplot [
color=blue,
solid,
mark=o,
mark options={solid}
]
table[row sep=crcr]{
1 0.42992843706819\\
2 0.236781091471918\\
3 0.128670455030046\\
4 0.0719661177485221\\
5 0.0421006325937349\\
6 0.0241921106872137\\
7 0.013781943093246\\
8 0.00762363762646751\\
9 0.00420728748893781\\
10 0.00229321995680742\\
11 0.00127205167649531\\
12 0.000713512972048673\\
13 0.000409319215427033\\
14 0.000232312382955809\\
15 0.000133454942130579\\
16 8.89019375912088e-05\\
17 7.20771147340466e-05\\
18 6.72420472934293e-05\\
19 6.57632833521202e-05\\
20 6.39438462104828e-05\\
21 6.32826937259935e-05\\
22 6.32415071654159e-05\\
23 6.32278108700502e-05\\
24 6.3232110469498e-05\\
25 6.32419802845572e-05\\
26 6.31586617124496e-05\\
27 6.31831217258071e-05\\
28 6.31872942657281e-05\\
29 6.32686394791312e-05\\
30 6.30885117435931e-05\\
};
\addlegendentry{$\widehat P_1$};

\addplot [
color=green!50!black,
solid,
mark=o,
mark options={solid}
]
table[row sep=crcr]{
1 0.043410621253962\\
2 0.00325161045776287\\
3 0.000189435995044932\\
4 1.27030391493213e-05\\
5 3.9404864803872e-06\\
6 3.85279112945228e-06\\
7 3.8526056563613e-06\\
8 3.85491464668788e-06\\
9 3.85289924751578e-06\\
10 3.85258143417745e-06\\
11 3.85301420009539e-06\\
12 3.85401722012563e-06\\
13 3.85340067675236e-06\\
14 3.85336759278756e-06\\
15 3.85386028341799e-06\\
16 3.85456692222014e-06\\
17 3.85450970735203e-06\\
18 3.85443847731832e-06\\
19 3.85457252238106e-06\\
20 3.8545175481955e-06\\
21 3.85491107058122e-06\\
22 3.85479959425301e-06\\
23 3.85474317741498e-06\\
24 3.85469770620749e-06\\
25 3.85471656147554e-06\\
26 3.8547660351178e-06\\
27 3.85467350767603e-06\\
28 3.85467603198196e-06\\
29 3.85550785401933e-06\\
30 3.85560815534587e-06\\
};
\addlegendentry{$\widehat P_5$};

\addplot [
color=red,
solid,
mark=o,
mark options={solid}
]
table[row sep=crcr]{
1 0.0189226398435413\\
2 0.000564255692704099\\
3 1.53226643443742e-05\\
4 4.04219183613066e-07\\
5 8.2577091393385e-08\\
6 8.19587404514969e-08\\
7 8.19904391063033e-08\\
8 8.18722139711352e-08\\
9 8.19937621799353e-08\\
10 8.17605491073142e-08\\
11 8.15278829351275e-08\\
12 8.16354807655249e-08\\
13 8.12222573586351e-08\\
14 8.10234707400678e-08\\
15 8.08579337255779e-08\\
16 8.06956585385658e-08\\
17 8.06533404993566e-08\\
18 8.0528806182052e-08\\
19 8.05929532126752e-08\\
20 8.06041659172464e-08\\
21 8.04697255355136e-08\\
22 8.04618739378176e-08\\
23 8.02831175801244e-08\\
24 8.00755529493216e-08\\
25 8.00603876391293e-08\\
26 8.00860081070553e-08\\
27 8.01036794001627e-08\\
28 8.00078981561397e-08\\
29 7.98433019781558e-08\\
30 7.99182810460797e-08\\
};
\addlegendentry{$\widehat P_{10}$};

\end{semilogyaxis}
\end{tikzpicture}%
  \caption{ Convergence of the Preconditioned GMRES using approximate
    iterates in $\cT_{5}(\cV)$ (top) or $\cT_{10}(\cV)$ (bottom), and
    using preconditioners $\widehat P_{r}$ (resp. $P_{r}$) constructed
    using ALG-P (resp. ALG-G) with different ranks $r \in \{1,5,10\}$.
    Comparison with the mean based preconditioner $P_E$.}
  \label{fig:largeleaf_conv_sol}
\end{figure}

In the following, we consider preconditioners with imposed sparsity
with $\gamma = 30\%$. Preconditioners constructed with ALG-P and ALG-G
are compared with mean-based preconditioner $P_{E}$ which is
classically used in the context of stochastic Galerkin methods.
$P_{E}$ is the sparse approximate inverse of the rank-one operator
associated with a mean-value of the parameters: $(\kappa,\eta) =
(5.5,600)$. A reference solution $u_{\text{ref}}$ is computed with the
  projection algorithm from \cite{Giraldi2012}, which yields a
  relative residual of $2.40~10^{-15}$ after 20 iterations. 
 In Figure \ref{fig:largeleaf_conv_sol}, we illustrate the convergence of the
  preconditioned GMRES for different preconditioners by plotting 
   the relative error $ \varepsilon(u^{(k)})$ between the $k$-th iterate $u^{(k)}$ of 
  GMRES and the reference solution $u_{\text{ref}}$, defined by
  \begin{equation}\label{eq:relative_err}
    \varepsilon(u^{(k)}) = \frac{\norm{u^{(k)} - u_{\text{ref}}}}{\norm{u_{\text{ref}}}}.
  \end{equation} 
  
\begin{figure}[h]
  \centering
  \setlength\figureheight{0.35\textwidth}
  \setlength\figurewidth{0.5\textwidth}
%
%
%
\begin{tikzpicture}

\begin{semilogyaxis}[%
width=\figurewidth,
height=\figureheight,
scale only axis,
xmin=0,
xmax=30,
xlabel={Iterations},
ymin=0.0001,
ymax=1,
yminorticks=true,
ylabel={$\widetilde \varepsilon(u^{(k)};I)$},
legend style={draw=black,fill=white,legend cell align=left}
]
\addplot [
color=blue,
solid,
mark=square,
mark options={solid}
]
table[row sep=crcr]{
1 0.32847631048525\\
2 0.134351986609456\\
3 0.0600299465138571\\
4 0.0289511023924653\\
5 0.0163669718506621\\
6 0.012514127787646\\
7 0.0108065852058826\\
8 0.0105819800752415\\
9 0.0105445077369007\\
10 0.0105024138236855\\
11 0.0105124405386724\\
12 0.0104843341981763\\
13 0.010447540844478\\
14 0.0104515955467427\\
15 0.0104508234122201\\
16 0.0104508107741224\\
17 0.0104529089190451\\
18 0.0104255954743443\\
19 0.0104220233039571\\
20 0.0104228979446605\\
21 0.0104184931752584\\
22 0.0104185603074629\\
23 0.0104185406066461\\
24 0.0104186606458462\\
25 0.0104194808635926\\
26 0.0104167844842653\\
27 0.0104155044984089\\
28 0.0104086018397824\\
29 0.0104035132424552\\
30 0.0104000536563979\\
};
\addlegendentry{$P_E$};

\addplot [
color=blue,
solid,
mark=+,
mark options={solid}
]
table[row sep=crcr]{
1 0.333650085773467\\
2 0.164315048536227\\
3 0.0819366433613296\\
4 0.0412997878979415\\
5 0.0227833468149929\\
6 0.0130387182270147\\
7 0.00847850725549274\\
8 0.00683803247464777\\
9 0.00643779838779507\\
10 0.00630885689945693\\
11 0.00623945175991754\\
12 0.00616137197064951\\
13 0.00612762233444384\\
14 0.00614566437399925\\
15 0.00615832873745128\\
16 0.00613472314709326\\
17 0.00613887876443177\\
18 0.00615453854567774\\
19 0.00616981384386193\\
20 0.00618025968102479\\
21 0.00615638002046867\\
22 0.0061310878675775\\
23 0.00613245143413026\\
24 0.00613566378885657\\
25 0.00614108509130885\\
26 0.00614215266958234\\
27 0.00614405293823059\\
28 0.00614494496793896\\
29 0.00614654901863171\\
30 0.00614504375333013\\
};
\addlegendentry{$P_1$};

\addplot [
color=green!50!black,
solid,
mark=+,
mark options={solid}
]
table[row sep=crcr]{
1 0.0743653519680228\\
2 0.0244960989880384\\
3 0.022479232650851\\
4 0.0222864407061029\\
5 0.0223197174006586\\
6 0.0223184468867023\\
7 0.0223159507725032\\
8 0.0223184653368282\\
9 0.0223187841594563\\
10 0.0223114251435034\\
11 0.0223122823248788\\
12 0.0223151807084879\\
13 0.0223133895090439\\
14 0.0223089013210963\\
15 0.0223056115658246\\
16 0.0223054561938952\\
17 0.0222942718086905\\
18 0.0222928758746492\\
19 0.0222818850170062\\
20 0.0222808250303462\\
21 0.0222741580412913\\
22 0.0222721561198781\\
23 0.0222722619103776\\
24 0.0222721831158125\\
25 0.0222722844726956\\
26 0.0222723121910356\\
27 0.0222721108514805\\
28 0.0222720510154841\\
29 0.0222732652285396\\
30 0.0222698542420441\\
};
\addlegendentry{$P_5$};

\addplot [
color=red,
solid,
mark=+,
mark options={solid}
]
table[row sep=crcr]{
1 0.075736161553269\\
2 0.0143152988504581\\
3 0.0122177001854527\\
4 0.0122050518232444\\
5 0.0122042166075011\\
6 0.0121734612882122\\
7 0.0121807490076399\\
8 0.0121627623607229\\
9 0.0121771586171787\\
10 0.0121778734962793\\
11 0.0121780602657833\\
12 0.0121715326501736\\
13 0.0121699024290708\\
14 0.0121641299651748\\
15 0.012164042138988\\
16 0.0121642167419366\\
17 0.0121652116663374\\
18 0.012165502744452\\
19 0.0121572364251628\\
20 0.0121386235141541\\
21 0.0121384428748949\\
22 0.0121383363748232\\
23 0.0121358133328509\\
24 0.0121323123767292\\
25 0.0121326154207896\\
26 0.012135424414044\\
27 0.0121360453649746\\
28 0.0121333924354804\\
29 0.0121325671137403\\
30 0.0121332185016783\\
};
\addlegendentry{$P_{10}$};

\addplot [
color=blue,
solid,
mark=o,
mark options={solid}
]
table[row sep=crcr]{
1 0.333650085773467\\
2 0.164315048536227\\
3 0.0819366433613345\\
4 0.0412997878979431\\
5 0.022783346815003\\
6 0.0130387496734606\\
7 0.0084789930342427\\
8 0.00683388804872603\\
9 0.00643100290281872\\
10 0.0063018203368752\\
11 0.00623277488353759\\
12 0.00615495398908383\\
13 0.00612135653987959\\
14 0.00613932392735899\\
15 0.00615179229567917\\
16 0.00612896205957837\\
17 0.00613293098415058\\
18 0.00614843145608281\\
19 0.00616390255248199\\
20 0.00617465315099533\\
21 0.00615722952538298\\
22 0.00612530502679276\\
23 0.00612606866533804\\
24 0.00612903816899436\\
25 0.00613433905665509\\
26 0.00613549590738785\\
27 0.00613725867161653\\
28 0.0061383050334128\\
29 0.00613981400847321\\
30 0.00613839044232127\\
};
\addlegendentry{$\widehat P_1$};

\addplot [
color=green!50!black,
solid,
mark=o,
mark options={solid}
]
table[row sep=crcr]{
1 0.0889271459627707\\
2 0.00731923082709633\\
3 0.000978039199131082\\
4 0.000892316357319528\\
5 0.000898690702071337\\
6 0.000896510875865551\\
7 0.000897565540212147\\
8 0.000897246306778589\\
9 0.000900790063429356\\
10 0.000898158743509488\\
11 0.000898788604146397\\
12 0.000899484218127232\\
13 0.000899915727277629\\
14 0.00089986195176326\\
15 0.000899245047024307\\
16 0.000899210990467709\\
17 0.000898902503581571\\
18 0.000898670846435317\\
19 0.00090236103260641\\
20 0.000897812714408924\\
21 0.000898346290034927\\
22 0.000898864307234367\\
23 0.00089893155565746\\
24 0.000898326610483469\\
25 0.000898823981897642\\
26 0.000899052303835529\\
27 0.00089910225885878\\
28 0.000899124625201386\\
29 0.000899509844153288\\
30 0.000899586316146805\\
};
\addlegendentry{$\widehat P_5$};

\addplot [
color=red,
solid,
mark=o,
mark options={solid}
]
table[row sep=crcr]{
1 0.0767262663560382\\
2 0.00255911005338633\\
3 0.000447266878087613\\
4 0.000445550197167832\\
5 0.000444895977372418\\
6 0.000436726791681985\\
7 0.000431974017024714\\
8 0.000429017717822035\\
9 0.000430791166348289\\
10 0.000427612695223306\\
11 0.000427930057050969\\
12 0.000427941081541325\\
13 0.000427284079116874\\
14 0.000427003730534062\\
15 0.000426860337364381\\
16 0.000426609540101152\\
17 0.000426728355172724\\
18 0.000426981311864664\\
19 0.000426952227509035\\
20 0.000426833722797544\\
21 0.00042684387879861\\
22 0.000426834880641779\\
23 0.000426868731638643\\
24 0.000426801273587294\\
25 0.000426786547360969\\
26 0.000426787861166632\\
27 0.000426763455406592\\
28 0.000426757652806923\\
29 0.000426751195984269\\
30 0.000426730380162246\\
};
\addlegendentry{$\widehat P_{10}$};

\end{semilogyaxis}
\end{tikzpicture}%
%
%
%
\begin{tikzpicture}

\begin{semilogyaxis}[%
width=\figurewidth,
height=\figureheight,
scale only axis,
xmin=0,
xmax=30,
xlabel={Iterations},
ymin=0.0001,
ymax=1,
yminorticks=true,
ylabel={$\varepsilon(u^{(k)})$},
legend style={draw=black,fill=white,legend cell align=left}
]
\addplot [
color=blue,
solid,
mark=square,
mark options={solid}
]
table[row sep=crcr]{
1 0.429158943473958\\
2 0.193924386555288\\
3 0.0881700553015909\\
4 0.0404095205676266\\
5 0.0191178928523393\\
6 0.0107669391559102\\
7 0.00752185987776282\\
8 0.0067955268152443\\
9 0.00662598662639562\\
10 0.00655272277776467\\
11 0.00655684465374879\\
12 0.00655630300154933\\
13 0.0065619491578679\\
14 0.00655953836048073\\
15 0.00656075779507579\\
16 0.00656073565122061\\
17 0.00655877837983758\\
18 0.00655693407536949\\
19 0.00655558513004954\\
20 0.00655546196484764\\
21 0.00655517403767528\\
22 0.00655517814663464\\
23 0.00655510429463314\\
24 0.00655517506253689\\
25 0.00655535779451559\\
26 0.00655429425936867\\
27 0.00655329365098749\\
28 0.00655058732142624\\
29 0.00654840847902228\\
30 0.00654643683192794\\
};
\addlegendentry{$P_E$};

\addplot [
color=blue,
solid,
mark=+,
mark options={solid}
]
table[row sep=crcr]{
1 0.42992843706819\\
2 0.236781091471918\\
3 0.128670509484902\\
4 0.0719687966549017\\
5 0.042131004380971\\
6 0.024491028603022\\
7 0.0148732302433656\\
8 0.0105147016028054\\
9 0.00917398726832181\\
10 0.00874480903476679\\
11 0.00856069625972374\\
12 0.00830798554908\\
13 0.00818421994861361\\
14 0.00816184273487837\\
15 0.00815864384970394\\
16 0.00816401187895246\\
17 0.00815976880554773\\
18 0.00815201773357234\\
19 0.00813037119186941\\
20 0.00812455141007925\\
21 0.00811170181003298\\
22 0.00808894256520563\\
23 0.00806997702683564\\
24 0.00802090934389001\\
25 0.00800975697202343\\
26 0.00800644304952495\\
27 0.00798158809492071\\
28 0.00798310713606205\\
29 0.00798150771356039\\
30 0.00798073966975526\\
};
\addlegendentry{$P_1$};

\addplot [
color=green!50!black,
solid,
mark=+,
mark options={solid}
]
table[row sep=crcr]{
1 0.0941624671358326\\
2 0.0142771640885178\\
3 0.00372001721389246\\
4 0.0033121753754568\\
5 0.00327970067801952\\
6 0.00328475303478438\\
7 0.0032832797499173\\
8 0.00328318106738159\\
9 0.00328324043748591\\
10 0.00328282637300127\\
11 0.0032832537571923\\
12 0.00328497582476427\\
13 0.00328492085746962\\
14 0.00328502664573372\\
15 0.00328577344950372\\
16 0.00328626756786768\\
17 0.00328536912228475\\
18 0.00328505457767686\\
19 0.00328496522913585\\
20 0.0032845363238357\\
21 0.0032841800522554\\
22 0.00328395465268379\\
23 0.00328418266068468\\
24 0.00328418656887807\\
25 0.00328414203166666\\
26 0.00328415471081687\\
27 0.00328402116457207\\
28 0.00328401615870062\\
29 0.00328369306635182\\
30 0.00328104041220835\\
};
\addlegendentry{$P_5$};

\addplot [
color=red,
solid,
mark=+,
mark options={solid}
]
table[row sep=crcr]{
1 0.0556955195383315\\
2 0.00655061148149525\\
3 0.00406226286321615\\
4 0.00404406773584782\\
5 0.00404108967391761\\
6 0.00404385409669075\\
7 0.0040428691167917\\
8 0.00404302743469598\\
9 0.0040425770107057\\
10 0.00404251396266575\\
11 0.0040424531555702\\
12 0.0040424351176822\\
13 0.00404112101761145\\
14 0.00403715747167926\\
15 0.00403715709119155\\
16 0.00403708590799767\\
17 0.00403709743753427\\
18 0.00403710736864559\\
19 0.00403706941535468\\
20 0.00403666313764792\\
21 0.0040365420766409\\
22 0.00403667409445501\\
23 0.00403673033173021\\
24 0.00403666486965566\\
25 0.00403675289578223\\
26 0.00403638389275186\\
27 0.00403642758422167\\
28 0.00403620929057795\\
29 0.00403612970978229\\
30 0.00403607185446815\\
};
\addlegendentry{$P_{10}$};

\addplot [
color=blue,
solid,
mark=o,
mark options={solid}
]
table[row sep=crcr]{
1 0.42992843706819\\
2 0.236781091471918\\
3 0.12867050948491\\
4 0.0719687966549024\\
5 0.0421310043809724\\
6 0.0244910337041273\\
7 0.0148731898108975\\
8 0.0105114575304236\\
9 0.00916870679004603\\
10 0.00873982097691359\\
11 0.00855671875534857\\
12 0.00830494628678802\\
13 0.00818142353525232\\
14 0.00815868809939376\\
15 0.00815536381089566\\
16 0.00816075985496901\\
17 0.00815645442689757\\
18 0.00814843116304996\\
19 0.00812652318697484\\
20 0.00812071355302719\\
21 0.00810856022891975\\
22 0.0080854102847881\\
23 0.00806609530727728\\
24 0.00801697408690945\\
25 0.00800587018238534\\
26 0.00800248837898055\\
27 0.00797738569143126\\
28 0.00797891651693714\\
29 0.00797729559171634\\
30 0.00797659504210993\\
};
\addlegendentry{$\widehat P_1$};

\addplot [
color=green!50!black,
solid,
mark=o,
mark options={solid}
]
table[row sep=crcr]{
1 0.043410621253962\\
2 0.00330019657469539\\
3 0.000723833573895744\\
4 0.000705244367550545\\
5 0.000704221134906199\\
6 0.000705242839717393\\
7 0.000705131810274106\\
8 0.000705152618076653\\
9 0.000705819758970892\\
10 0.000704959465422618\\
11 0.000704805427053817\\
12 0.000704802208786707\\
13 0.000704775513261455\\
14 0.000704752192235011\\
15 0.000704666615769654\\
16 0.000704607600163609\\
17 0.00070458388693336\\
18 0.000704302329437687\\
19 0.000704946131834381\\
20 0.000704093696626489\\
21 0.000704099858543493\\
22 0.000704191466609824\\
23 0.000704190110123061\\
24 0.000704192003631242\\
25 0.000704204709636677\\
26 0.000704238572511431\\
27 0.000704242349805227\\
28 0.000704244405791343\\
29 0.000704264397982547\\
30 0.000704218987291762\\
};
\addlegendentry{$\widehat P_5$};

\addplot [
color=red,
solid,
mark=o,
mark options={solid}
]
table[row sep=crcr]{
1 0.0189227860236283\\
2 0.000585883642294805\\
3 0.000138862933447314\\
4 0.000137995477959175\\
5 0.000136289279131353\\
6 0.000138383346891256\\
7 0.000138045544062411\\
8 0.000138169504448673\\
9 0.000138368530773669\\
10 0.000139635901674672\\
11 0.000139799285053435\\
12 0.000139812492137298\\
13 0.000139468253682434\\
14 0.000139499579599182\\
15 0.000139505933327069\\
16 0.000139526733404065\\
17 0.00013952216339755\\
18 0.000139615836784086\\
19 0.000139623051548641\\
20 0.000139669956556833\\
21 0.00013967006022306\\
22 0.000139673551438144\\
23 0.000139692198289339\\
24 0.000139753839882805\\
25 0.000139781236810756\\
26 0.000139782338004014\\
27 0.000139755988357834\\
28 0.000139746068819785\\
29 0.000139771515158977\\
30 0.000139779579458205\\
};
\addlegendentry{$\widehat P_{10}$};

\end{semilogyaxis}
\end{tikzpicture}%
%
%
%
\begin{tikzpicture}

\begin{semilogyaxis}[%
width=\figurewidth,
height=\figureheight,
scale only axis,
xmin=0,
xmax=30,
xlabel={Iterations},
ymin=0.0001,
ymax=1,
yminorticks=true,
ylabel={$\widetilde \varepsilon(u^{(k)};P)$},
legend style={draw=black,fill=white,legend cell align=left}
]
\addplot [
color=blue,
solid,
mark=square,
mark options={solid}
]
table[row sep=crcr]{
1 0.307423517480681\\
2 0.118654225910475\\
3 0.0495498822106539\\
4 0.0216660386200253\\
5 0.0101342389398768\\
6 0.00584717936146093\\
7 0.00464068451189854\\
8 0.00443633279613137\\
9 0.00439012308231164\\
10 0.00435441949469978\\
11 0.00435739366363809\\
12 0.00435414519540056\\
13 0.00434028239083569\\
14 0.004343018040242\\
15 0.00434300019124039\\
16 0.0043430882942618\\
17 0.00434000130101132\\
18 0.00433336860916073\\
19 0.00433183769072235\\
20 0.00433277834380718\\
21 0.0043327329908239\\
22 0.00433276191195882\\
23 0.00433275291789583\\
24 0.00433281352069596\\
25 0.00433272874259781\\
26 0.00433246770102674\\
27 0.00433223061177542\\
28 0.00433352721028779\\
29 0.00433375059240318\\
30 0.00433444959240235\\
};
\addlegendentry{$P_E$};

\addplot [
color=blue,
solid,
mark=+,
mark options={solid}
]
table[row sep=crcr]{
1 0.3771731818437\\
2 0.166464150456228\\
3 0.0828777432044831\\
4 0.0405842925650412\\
5 0.0213907056048375\\
6 0.0120470522806327\\
7 0.00645013712056718\\
8 0.00500705023205188\\
9 0.00474238066741502\\
10 0.00478701578162379\\
11 0.00474929503535454\\
12 0.0046907121947947\\
13 0.00469052725305355\\
14 0.00473500254283332\\
15 0.00475660538707674\\
16 0.00471425960996673\\
17 0.00471743107132614\\
18 0.00473446046421701\\
19 0.00474777581244238\\
20 0.00476276727087379\\
21 0.00472637570603394\\
22 0.00469073355679443\\
23 0.00469024386651789\\
24 0.00467090500108941\\
25 0.00467673536936395\\
26 0.0046788301390854\\
27 0.00467107176826118\\
28 0.00467210487020546\\
29 0.00466739388574513\\
30 0.004665855689014\\
};
\addlegendentry{$P_1$};

\addplot [
color=green!50!black,
solid,
mark=+,
mark options={solid}
]
table[row sep=crcr]{
1 0.07963639762115\\
2 0.0117588322604053\\
3 0.00361132641075731\\
4 0.00337922805040173\\
5 0.00335408407793033\\
6 0.00335596650341889\\
7 0.00335575369059678\\
8 0.00335579003186791\\
9 0.00335585199629424\\
10 0.00335534877245383\\
11 0.00335566204813191\\
12 0.00335691102703887\\
13 0.00335698643123645\\
14 0.00335680370828309\\
15 0.00335743453442526\\
16 0.00335768861648097\\
17 0.00335704939164015\\
18 0.00335668383671086\\
19 0.0033563245614183\\
20 0.00335598342156267\\
21 0.0033555891965745\\
22 0.0033553636800465\\
23 0.00335558534193608\\
24 0.00335559618031441\\
25 0.00335555371094216\\
26 0.00335553962943031\\
27 0.00335541546332847\\
28 0.00335540017180417\\
29 0.00335470704900174\\
30 0.00335224363888636\\
};
\addlegendentry{$P_5$};

\addplot [
color=red,
solid,
mark=+,
mark options={solid}
]
table[row sep=crcr]{
1 0.0502754020812139\\
2 0.00600896726198259\\
3 0.00405758642798417\\
4 0.00404524569365444\\
5 0.00404369761858536\\
6 0.00404520301030651\\
7 0.00404463202717378\\
8 0.00404508734594017\\
9 0.00404441871661899\\
10 0.00404435112247382\\
11 0.00404428008550836\\
12 0.0040444751353026\\
13 0.00404301298454281\\
14 0.00403898615150742\\
15 0.0040389986706684\\
16 0.0040389456098228\\
17 0.00403893851118208\\
18 0.00403894017750567\\
19 0.004039027168357\\
20 0.00403867266609273\\
21 0.00403854952093615\\
22 0.00403866385926766\\
23 0.00403872181964627\\
24 0.00403872926429699\\
25 0.00403875357230333\\
26 0.00403848289859435\\
27 0.00403851295807301\\
28 0.00403833221301772\\
29 0.00403827744459436\\
30 0.00403825021488463\\
};
\addlegendentry{$P_{10}$};

\addplot [
color=blue,
solid,
mark=o,
mark options={solid}
]
table[row sep=crcr]{
1 0.3771731818437\\
2 0.166464150456228\\
3 0.0828777432044829\\
4 0.0405842925650416\\
5 0.0213907056048392\\
6 0.0120471491222011\\
7 0.00645029633859103\\
8 0.00500576256241915\\
9 0.00473972409003457\\
10 0.00478427982566549\\
11 0.00474681501436219\\
12 0.00468854032366388\\
13 0.00468836193790016\\
14 0.0047327701774486\\
15 0.0047545178681283\\
16 0.00471291322329917\\
17 0.00471618282183168\\
18 0.00473285980641229\\
19 0.00474583261989078\\
20 0.00476183884658615\\
21 0.00473407956707565\\
22 0.00468904373068218\\
23 0.00468804255161421\\
24 0.00466841961062866\\
25 0.0046742278797167\\
26 0.00467637220334647\\
27 0.00466840247431969\\
28 0.00466943774928952\\
29 0.0046647311101418\\
30 0.00466327305872599\\
};
\addlegendentry{$\widehat P_1$};

\addplot [
color=green!50!black,
solid,
mark=o,
mark options={solid}
]
table[row sep=crcr]{
1 0.0403902140617849\\
2 0.00295780524181447\\
3 0.000724916988342065\\
4 0.000708154792272189\\
5 0.000707049277313508\\
6 0.000707921351248478\\
7 0.000707859227021839\\
8 0.0007079267929118\\
9 0.000708626162085893\\
10 0.00070773989312769\\
11 0.000707600346513676\\
12 0.000707603311160372\\
13 0.000707575490135669\\
14 0.000707546374300698\\
15 0.000707462372999757\\
16 0.000707406389422424\\
17 0.000707384306948631\\
18 0.000707099611163453\\
19 0.000707751881679149\\
20 0.000706898784061724\\
21 0.000706897885902679\\
22 0.000706974907040863\\
23 0.000706973676743915\\
24 0.000706976003951047\\
25 0.000706993182955157\\
26 0.000707025576539275\\
27 0.000707031245159361\\
28 0.000707033310924293\\
29 0.000707049772051896\\
30 0.000706997806968372\\
};
\addlegendentry{$\widehat P_5$};

\addplot [
color=red,
solid,
mark=o,
mark options={solid}
]
table[row sep=crcr]{
1 0.0186765183075547\\
2 0.000569056954310388\\
3 0.000137830600286533\\
4 0.000137029992116387\\
5 0.000135544164730749\\
6 0.000137527248756076\\
7 0.000137307198468682\\
8 0.000137402085034885\\
9 0.000137594016646755\\
10 0.000138767127201141\\
11 0.000138924171418697\\
12 0.000138937042066615\\
13 0.00013860327781265\\
14 0.000138632399143645\\
15 0.000138646752633956\\
16 0.000138668297129235\\
17 0.000138665686001683\\
18 0.000138755883344\\
19 0.000138761813803659\\
20 0.000138804964033728\\
21 0.000138804869081216\\
22 0.000138808020777598\\
23 0.000138827419271919\\
24 0.000138883958879101\\
25 0.000138909773298336\\
26 0.000138910795207457\\
27 0.000138883857236966\\
28 0.000138874366869087\\
29 0.000138899765742573\\
30 0.000138907117073663\\
};
\addlegendentry{$\widehat P_{10}$};

\end{semilogyaxis}
\end{tikzpicture}%
  \caption{ Convergence of the Preconditioned GMRES using approximate
    iterates in $\cT_{5}(\cV)$ for different error estimators: 
    relative residual norm (top), relative error in solution norm (middle) and 
    preconditioned relative residual norm (bottom). Use of preconditioners $\widehat P_{r}$ (resp. $P_{r}$) constructed
    with ALG-P (resp. ALG-G) with different ranks $r \in \{1,5,10\}$.
    Comparison with the mean based preconditioner $P_E$.}
  \label{fig:largeleaf_conv_sol_err_estim}
\end{figure}
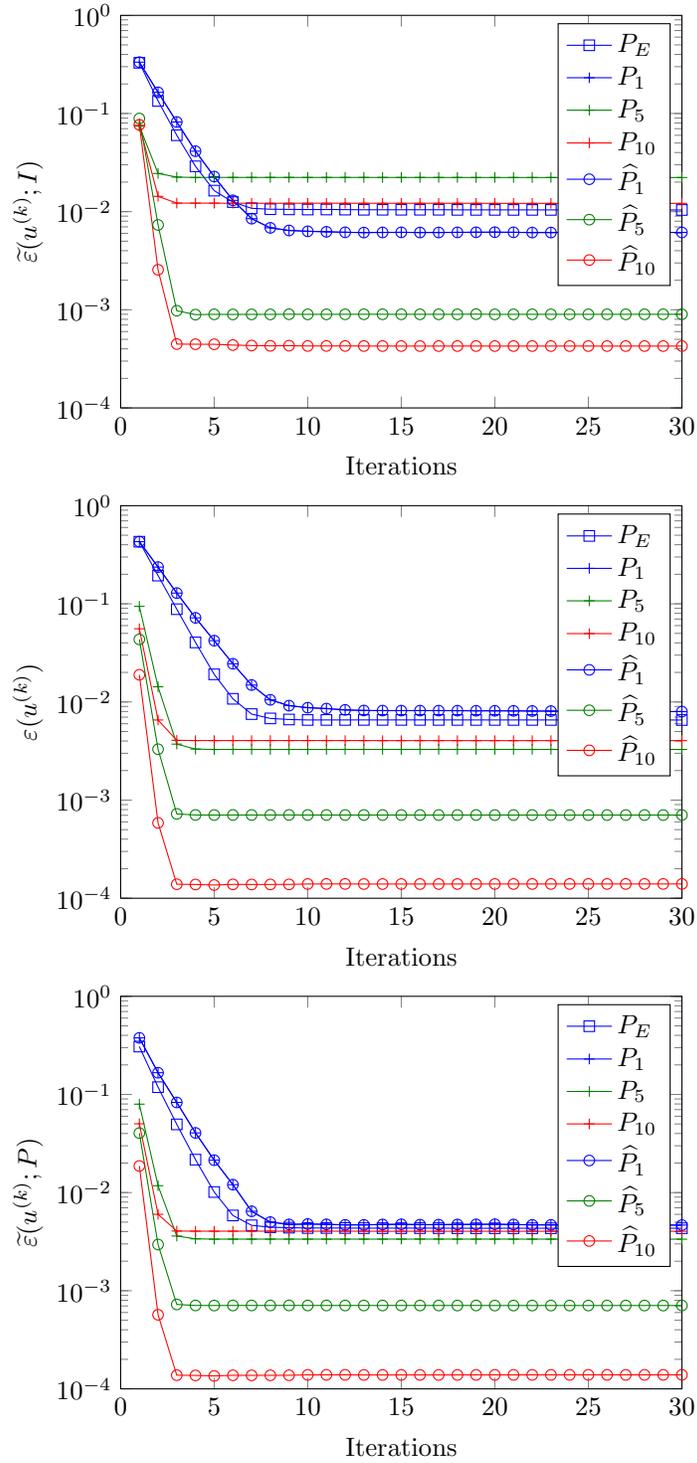

  We observe that all the rank-one preconditioners give similar
  convergences, and that ALG-P greatly improves the convergence rate of the
  preconditioned GMRES algorithm. 
  We observe a stagnation of the relative error at a certain
  precision. This precision depends on the tensor subset in which the
  iterates are approximated (here $\cT_{10}(\cV)$ or $\cT_{5}(\cV)$)
  and also on the conditioning of the operator. 
  The final precision
  can first be improved by introducing larger tensor subsets for the
  approximation of iterates (compare the precisions obtained using
  $\cT_{10}(\cV)$ or $\cT_{5}(\cV)$). This was already observed in
  \cite{kressner_low-rank_2011,Ballani2013}. Also, the final precision
  can be improved by using a better preconditioner. Figures
  \ref{fig:largeleaf_conv_sol} (top and bottom) illustrate that the
  relative error for an approximation in $\cT_5(\cV)$ and a
  preconditioner $\widehat P_{10}$ obtained with ALG-P is
  $1.4~10^{-4}$, while the relative error for an approximation in
  $\cT_{10}(\cV)$ and a preconditioner $P_{1}$ is $2.3~10^{-4}$.
  Concerning the case where the iterates are approximated in
  $\cT_{10}(\cV)$ (Figure \ref{fig:largeleaf_conv_sol} bottom), we
  observe that the final relative error slightly increases with the
  rank $r$ of the preconditioner $P_r$ constructed by ALG-G. This
  reflects the fact that when using ALG-G, increasing $r$ only
  slightly improves the quality of the preconditioner and may result
  in a deterioration of the final precision measured in solution norm.
  However,
  the convergence rate of GMRES clearly increases with
  $r$.
In practice, when a preconditioner $P$ is available, the error can be estimated by computing 
 the norm of the preconditioned relative
  residual defined by
  \begin{equation*}
    \widetilde \varepsilon(u^{(k)};P) = \frac{\norm{P(b-Au^{(k)})}}{\norm{Pb}}. 
  \end{equation*}
For good preconditioners,  $\widetilde \varepsilon(u^{(k)};P)$ may be a good estimate of the relative error $\varepsilon(u^{(k)})$ in solution norm. 
  This is illustrated on Figure
  \ref{fig:largeleaf_conv_sol_err_estim} where we can see that
  $\widetilde \varepsilon(u^{(k)};P)$ gives almost the same error as
  $\varepsilon(u^{(k)})$ when the preconditioner $P$ is a good
  approximation of $A^{-1}$ (e.g. $P = P_5$, $P_{10}$, $\widehat P_5$, $\widehat
  P_{10}$). 

\paragraph{Influence of sparsity}
\label{sec:influence-gamma}

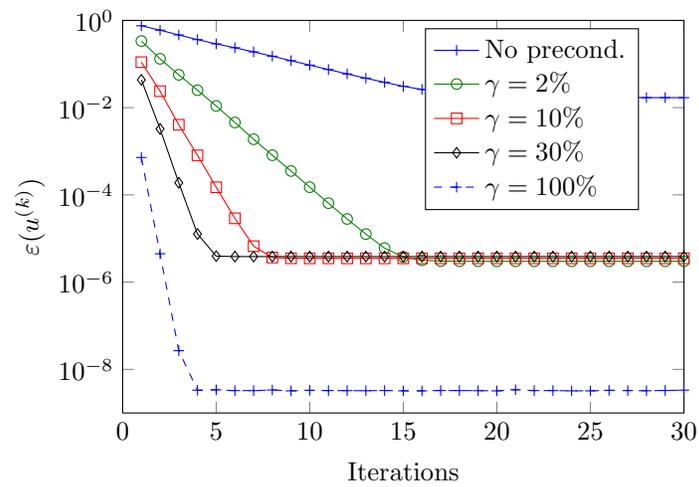
\begin{figure}[h]
  \centering
  \setlength\figureheight{0.35\textwidth}
  \setlength\figurewidth{0.5\textwidth}
%
%
%
%
\begin{tikzpicture}

\begin{semilogyaxis}[%
width=\figurewidth,
height=\figureheight,
scale only axis,
xmin=0,
xmax=30,
xlabel={Iterations},
ymin=1e-09,
ymax=1,
yminorticks=true,
ylabel={$\varepsilon(u^{(k)})$},
legend style={at={(0.92,0.98)},draw=black,fill=white,legend cell align=left}
]
\addplot [
color=blue,
solid,
mark=+,
mark options={solid}
]
table[row sep=crcr]{
1 0.752458666511062\\
2 0.596651548132291\\
3 0.463522154396855\\
4 0.362319851161486\\
5 0.289735477963995\\
6 0.236203733829995\\
7 0.190078438202346\\
8 0.151128698173056\\
9 0.11939052254939\\
10 0.0943566054954811\\
11 0.0748594455836435\\
12 0.0593275594351322\\
13 0.0472752380463031\\
14 0.0380905982667648\\
15 0.0308614146774255\\
16 0.0261284510048281\\
17 0.022580297512802\\
18 0.0199564920129199\\
19 0.0184254199365102\\
20 0.0177712821767183\\
21 0.0173058042925715\\
22 0.0170941990257486\\
23 0.0170855853829567\\
24 0.0170598300802174\\
25 0.0170427995671194\\
26 0.0170213926783831\\
27 0.0169909060649116\\
28 0.0169657108708949\\
29 0.0169509916764406\\
30 0.0169492828148351\\
};
\addlegendentry{No precond.};

\addplot [
color=green!50!black,
solid,
mark=o,
mark options={solid}
]
table[row sep=crcr]{
1 0.336197136889593\\
2 0.13167434795156\\
3 0.0565734105735576\\
4 0.0250939806185684\\
5 0.0109353765943859\\
6 0.00456513254011281\\
7 0.00188280228884709\\
8 0.000810395352126042\\
9 0.000352046767435025\\
10 0.000149828676523594\\
11 6.45073392499918e-05\\
12 2.80167638281516e-05\\
13 1.26208530433342e-05\\
14 6.02398125912157e-06\\
15 3.77836773473066e-06\\
16 3.17211681906938e-06\\
17 3.03564044192676e-06\\
18 3.0068058927854e-06\\
19 2.99928272290413e-06\\
20 2.99927474856554e-06\\
21 2.99926814076408e-06\\
22 2.99980031794422e-06\\
23 2.99976446653333e-06\\
24 2.99911839285193e-06\\
25 2.99930533075556e-06\\
26 2.99908034897152e-06\\
27 2.9991557673352e-06\\
28 2.99919363378566e-06\\
29 2.99958687534687e-06\\
30 2.9993059465757e-06\\
};
\addlegendentry{$\gamma = 2\%$};

\addplot [
color=red,
solid,
mark=square,
mark options={solid}
]
table[row sep=crcr]{
1 0.110273350834182\\
2 0.023875394995532\\
3 0.00404945169040982\\
4 0.000806357009057141\\
5 0.000149492532931326\\
6 2.8998617879211e-05\\
7 6.65598556140714e-06\\
8 3.65608963895756e-06\\
9 3.51611831129894e-06\\
10 3.50640561509306e-06\\
11 3.50663985333213e-06\\
12 3.50544148857316e-06\\
13 3.50494924749468e-06\\
14 3.50464313479347e-06\\
15 3.50422929203411e-06\\
16 3.50492912736863e-06\\
17 3.5042156304346e-06\\
18 3.50512098095868e-06\\
19 3.50474650825625e-06\\
20 3.50458478248644e-06\\
21 3.50482065466456e-06\\
22 3.50466081133816e-06\\
23 3.50493282440245e-06\\
24 3.50546719885965e-06\\
25 3.5067920766076e-06\\
26 3.50607576413663e-06\\
27 3.50399214136352e-06\\
28 3.50395465047386e-06\\
29 3.50410410625572e-06\\
30 3.5061186366007e-06\\
};
\addlegendentry{$\gamma = 10\%$};

\addplot [
color=black,
solid,
mark=diamond,
mark options={solid}
]
table[row sep=crcr]{
1 0.0434106212539619\\
2 0.00325161045776105\\
3 0.000189435995059486\\
4 1.27030384864145e-05\\
5 3.94048441349896e-06\\
6 3.85279166871679e-06\\
7 3.85260650386587e-06\\
8 3.85491479739347e-06\\
9 3.85289841882941e-06\\
10 3.85258108619095e-06\\
11 3.85300805244894e-06\\
12 3.854016860206e-06\\
13 3.85339990611283e-06\\
14 3.85336776501976e-06\\
15 3.85386114810571e-06\\
16 3.85456799332875e-06\\
17 3.85451085090378e-06\\
18 3.85443840094575e-06\\
19 3.85457111711862e-06\\
20 3.85451727222316e-06\\
21 3.85490798438723e-06\\
22 3.85480063561023e-06\\
23 3.85474267120356e-06\\
24 3.85469818098586e-06\\
25 3.85471759974136e-06\\
26 3.85476498734532e-06\\
27 3.85467297679547e-06\\
28 3.85467535501632e-06\\
29 3.8555077856113e-06\\
30 3.85560696249266e-06\\
};
\addlegendentry{$\gamma = 30\%$};

\addplot [
color=blue,
dashed,
mark=+,
mark options={solid}
]
table[row sep=crcr]{
1 0.000720266519822012\\
2 4.46962620333699e-06\\
3 2.67884011529185e-08\\
4 3.32759704150292e-09\\
5 3.39093245229049e-09\\
6 3.26040900728922e-09\\
7 3.25347894019072e-09\\
8 3.36883341372934e-09\\
9 3.21082424615271e-09\\
10 3.31047065882032e-09\\
11 3.26345665893424e-09\\
12 3.25288655675377e-09\\
13 3.23783994873341e-09\\
14 3.27261899243982e-09\\
15 3.21525450109855e-09\\
16 3.21121898263454e-09\\
17 3.27152832172968e-09\\
18 3.26125106933677e-09\\
19 3.24477468294117e-09\\
20 3.21269068067794e-09\\
21 3.41292704999768e-09\\
22 3.26059286288813e-09\\
23 3.26231670314943e-09\\
24 3.21509868731631e-09\\
25 3.21411830276486e-09\\
26 3.30664220969821e-09\\
27 3.26828307882792e-09\\
28 3.23262339107452e-09\\
29 3.27792410265925e-09\\
30 3.33992930840636e-09\\
};
\addlegendentry{$\gamma = 100\%$};

\end{semilogyaxis}
\end{tikzpicture}%
  \caption{Convergence of Preconditioned GMRES using approximate
    iterates in $\cT_{10}(\cV)$ and using preconditioner $\widehat P_{5}$
    constructed using ALG-P with different sparsity levels: $\gamma =
    2,10,30,100 \%$.}
  \label{fig:influence-gamma}
\end{figure}

We consider the preconditioner $\widehat P_{5}$ obtained with ALG-P.
The convergence of the preconditioned GMRES solver using
approximations of iterates in $\cT_{10}(\cV)$ for different levels $\gamma$ of
sparsity is plotted in Figure \ref{fig:influence-gamma}.  We observe that the preconditioning greatly improves the convergence rate of GMRES and also the accuracy of the resulting approximation. 
  As we could have expected, increasing $\gamma$ improves the
preconditioner and therefore improves the convergence rate and the
quality of the resulting approximation. When $\gamma = 100\%$, that
means without imposed sparsity, the algorithm converges very fast and
the error $\varepsilon(u^{(k)})$ stagnates at a very low value of
$3~10^{-9}$ after only 4 iterations. When $\gamma = 2\%$, $10\%$ or
$30\%$, we observe that the error $\varepsilon(u^{(k)})$ stagnates at
about the same value $3~10^{-6}$, which is greater than for
$\gamma=100\%$ but also significantly lower than the final error of
$1.7~10^{-2}$ obtained when no preconditioner is used.

\clearpage
\section{Conclusion}

An algorithm has been proposed for the progressive construction of
low-rank approximations of the inverse of an operator given in
low-rank tensor format. This construction is based on an updated
greedy algorithm which consists in constructing a sequence of tensor
subspaces from successive rank-one corrections and in computing
projections (or approximate projections) in these tensor subspaces,
thus resulting in approximations in low-rank Tucker (or Hierarchical
Tucker) format. Some desired properties can be imposed on the
approximate inverse during the correction step, such as symmetry
(requiring the solution of Sylvester equations) or sparsity. For the
latter case, the algorithm relies on a straightforward adaptation of
the SParse Approximate Inverse method with adaptive selection of the
pattern. Compared to a direct approximation in low-rank tensor
subsets, the updated greedy algorithm has the advantages of being
adaptive and of allowing the reduction of the complexity of the
construction. Also, the projection step significantly improves the
quality of preconditioners which may be obtained with pure greedy
rank-one algorithms. 
Numerical examples have illustrated the ability
of the algorithm to provide good preconditioners for linear systems of
equations with a significant improvement of the convergence
properties of iterative solvers. The final precision which can be
obtained in fixed low-rank tensor subsets is enhanced as well.

Further investigations are needed in order to measure the quality of
the approximate inverse as a preconditioner of a linear system of
equations. This measure could provide pertinent stopping criteria for
the updated greedy algorithm, with a necessary balance between the
quality of the preconditioner and the computational complexity
(complexity of its construction and of algorithms for solving the
preconditioned system).

\bibliographystyle{siam}

\end{document}